\documentclass[11 pt, psamsfonts]{amsart}
\usepackage{amsmath}
\usepackage{amsthm}
\usepackage{amsfonts}
\usepackage{amssymb}
\usepackage{amscd,epsf, latexsym}
\usepackage[dvips]{graphicx}

\usepackage{eucal}
\newcommand{\1}{{1 \hspace{-0.35em} {\rm 1}}}
\newcommand{\I}{\mathcal{I}}
\newcommand{\sgn}{\ve}

\newcommand{\M}{\mathcal{M}}
\newcommand{\OO}{\mathcal{O}}
\newcommand{\Vn}{V^{\otimes n}}
\newcommand{\B}{\mathcal{B}}
\newcommand{\V}{\mathcal{V}}

\newcommand{\lan}{\langle}
\newcommand{\ra}{\rangle}
\newcommand{\End}{{\rm End}}
\newcommand{\Hom}{{\rm Hom}}
\newcommand{\ve}{\varepsilon}
\newcommand{\T}{\mathcal{T}}
\newcommand{\tq}{\tilde{q}}
\newcommand{\F}{\mathcal{F}}
\newcommand{\mF}{\mathfrak{F}}
\newcommand{\mV}{\mathfrak{V}}
\newcommand{\Dim}{{\rm Dim}}
\newcommand{\zero}{{\mathbf 0}}
\newcommand{\ga}{\gamma}
\newcommand{\Ga}{\Gamma}
\newcommand{\one}{{\rm Id}}

\newcommand{\Z}{\mathbb{Z}}
\newcommand{\g}{\mathfrak{g}}
\newcommand{\N}{\mathbb{N}}
\newcommand{\Q}{Q}
\newcommand{\C}{\mathbb{C}}
\newcommand{\R}{\mathbb{R}}

\newcommand{\la}{{\lambda}}
\newcommand{\La}{{\Lambda}}
\newcommand{\al}{\alpha}

\theoremstyle{plain}
\newtheorem{theorem}{Theorem}[section]
\newtheorem{lemma}[theorem]{Lemma}

\newtheorem{scholium}[theorem]{Scholium}
\newtheorem{prop}[theorem]{Proposition}
\newtheorem{cor}[theorem]{Corollary}
\theoremstyle{definition}
\newtheorem{definition}[theorem]{Definition}
\theoremstyle{remark}

\newtheorem{rmk}[theorem]{Remark}
\newtheorem{example}[theorem]{Example}

\calclayout

\title{On a Family of Non-Unitarizable Ribbon Categories}
\author{Eric C. Rowell}
\thanks{\textit{Department of Mathematics, Indiana University, Bloomington, IN 47405 USA}\\
email: \texttt{errowell@indiana.edu}}
\begin{document}

\begin{abstract}
We consider several families of categories.  The first are
quotients of Andersen's tilting module categories for quantum
groups of Lie type $B$ at odd roots of unity.  The second consists
of categories of type $BC$ constructed from idempotents in
$BMW$-algebras. Our main result is to show that these families
coincide as braided tensor categories using a recent theorem of
Tuba and Wenzl. By appealing to similar results of Blanchet and
Beliakova we obtain another interesting equivalence with these two
families of categories and the quantum group categories of Lie
type $C$ at odd roots of unity.  The morphism spaces in these
categories can be equipped with a Hermitian form, and we are able
to show that these categories are never unitary, and no braided
tensor category sharing the Grothendieck semiring common to these
families is unitarizable.

\end{abstract}

\maketitle

\section{Introduction}
The purpose of this paper is two-fold: to solve an open problem
regarding the unitarity of Hermitian ribbon categories arising
from quantum groups, and to make progress towards the
classification of finite ribbon categories.

To any simple Lie algebra $\g$ and a parameter $q$ with $q^2$ a
primitive $\ell$th root of unity one may associate a finite
semisimple Hermitian ribbon category $\F$ derived from
representations of quantum groups. A further property that $\F$
may have is \emph{unitarity}, which depends on the algebra $\g$
and the specific choice of $q$. In 1998 Wenzl \cite{wenzl} showed
that for $\g$ of simply-laced type there is always a choice $q$
that yields a unitary category, and for non-simply-laced types as
long as $\ell$ is divisible by 2 (resp. 3) for types $B$, $C$ and
$F$ (resp. $G$).  It was hoped that these divisibility conditions
could be removed by making a clever choice of $q$ or changing the
braiding, but whether this was possible remained a dark mystery.
This was the original motivation for this paper--to explore
unitarity for this family of type $B$, odd $\ell$ categories.

Among the other constructions of ribbon categories that are
currently known, one of the most interesting blends ideas from
operator algebras and link invariants and is essentially due to
Turaev and Wenzl \cite{TuraevWenzl2}.  Recently Tuba and Wenzl
\cite{TuW2} studied these families of categories and were able to
get a partial classification--determining the possible braiding
and monoidal structures from the Grothendieck semiring. We use
their result to identify the aforementioned family of Lie type
$B$, odd $\ell$ quantum group categories with certain Turaev-Wenzl
categories of ortho-symplectic $BC$ type at the level of braided
tensor categories.  Similar results were obtained by Beliakova
and Blanchet in \cite{BB}.  The main equivalence we establish
is just an extension of an equivalence Beliakova and Blanchet
observed to spin modules.  Combining their results with ours we
get as a corollary a rank-level type duality between the Lie type
$B$ and $C$ quantum group categories at odd roots of unity (see
Corollary \ref{ranklevel}).

By a thorough (but elementary) analysis of characters of the
Grothendieck semirings of these categories, we are able to show
that no Hermitian ribbon category with the same tensor product
rules as these categories can be unitary.  Thus we answer the
original question of unitarity for both quantum groups of Lie
types $B$ and $C$ at odd roots of unity.

The organization of this paper is as follows.  In Section 2 we
define the categorical terms of the subject and mention a few
results germane to the discussion. In Sections 3 and 4 we describe
the structure of the family of quantum group categories we are
concerned with and analyze the Grothendieck semiring and
characters.  This sets the stage Section 5 in which we consider
the representations of the braid group on morphism spaces and the
second family of categories we consider.  In Section 6 we
establish the equivalence between these two families of
categories.  In Section 7 we apply this equivalence to prove the
failure of unitarity.

\subsection*{Acknowledgements} The author would like to thank H. Wenzl, I. Tuba, N.
Wallach and Z. Wang for many useful discussions.  The author is
also grateful to a referee who pointed out a misinterpretation of
the multiplicity formula of \cite{AndPar}, and to Andersen and
S. Sawin for their help in fixing this error.  It recently came to
the author's attention that the (almost) character-preserving
involution of the Weyl alcove described below appeared in a
different (more general) form in papers by S. Sawin (\cite{Sawin})
and T. Le and V. Turaev (\cite{LeTuraev}) about the same time as
it was discovered by the author (see \cite{thesis}).  The
techniques used here are somewhat more elementary than those found
in these works.
\section{Ribbon Categories}

\subsection{Axioms}\label{cat}
In this subsection we outline the relevant categorical axioms.
We follow the paper \cite{TuraevWenzl2}, and refer to that paper or the books
by Turaev \cite{Tur} or Kassel \cite{K} for a complete treatment.\\

Let $\OO$ be a category defined over a subfield $k\subset\C$.
The following axioms are satisfied by a semisimple
Hermitian ribbon category.
\begin{enumerate}
\item[1.] A \textbf{monoidal category} has a tensor product $\otimes$
and an identity object $\1$ satisfying the triangle and pentagon axioms.
These guarantee that the tensor product is associative (at least up to
isomorphism) and that
\begin{equation*}\label{idob}
\1\otimes X\cong X\otimes\1\cong X
\end{equation*}
 for any object $X$.  We usually assume our categories are \emph{strict},
that is, that the associativity isomorphisms and the isomorphisms above are
the identity.
\item[2.] A category is \textbf{rigid} if there is a dual
module $X^*$ for each object $X$
and morphisms
$$b_X: \1 \rightarrow X\otimes X^*, d_X: X^*\otimes X \rightarrow \1$$
satisfying
\begin{eqnarray}
(\one_X\otimes d_X)(b_X\otimes \one_X)&=&\one_X \\
(d_X\otimes \one_{X^*})(\one_{X*}\otimes b_X)&=&\one_{X^*}.
\end{eqnarray}
\item[3.] An \textbf{Ab-category} is one in which all morphism spaces
are $\C$-vector spaces and the composition and tensor product of
morphisms are bilinear.
\item[4.] A \textbf{semisimple} category has the property that every object
$X$ is isomorphic to a finite direct sum of \emph{simple}
objects--that is, objects $X_i$ with $\End(X_i)\cong\C$--and that
the simple objects satisfy Schur's Lemma: $\dim\Hom(X_i,X_j)\in
\{0,1\}$.  $\OO$ is called \textbf{finite} if there are finitely
many isomorphism classes of simple objects.
\item[5.] A
\textbf{braiding} is a family of isomorphisms
$$c_{X,Y}: X\otimes Y \rightarrow Y\otimes X$$
satisfying
\begin{eqnarray}
c_{X,Y\otimes Z}&=&(\one_Y\otimes c_{X,Z})(c_{X,Y}\otimes \one_Z) \\
c_{X\otimes Y,Z}&=&(c_{X,Z}\otimes \one_Y)(\one_X\otimes c_{Y,Z})
\end{eqnarray}
\item[6.] A \textbf{twist} consists of isomorphisms
$$\theta_X: X \rightarrow X$$
To be compatible with the braiding and duality we must have:
\begin{eqnarray}
\theta_{X\otimes Y}&=&c_{Y,X}c_{X,Y}(\theta_X\otimes \theta_Y) \\
\theta_{X^*}&=&(\theta_X)^*
\end{eqnarray}
A rigid category is called \emph{balanced} if it has a twist.
\item[7.] A \textbf{Hermitian} category has a conjugation:
$$\dag: \Hom(X,Y)\rightarrow \Hom(Y,X)$$
such that $(f^\dag)^\dag=f$, $(f\otimes g)^\dag=f^\dag\otimes g^\dag$
and $(f\circ g)^\dag=g^\dag\circ f^\dag$.  On $\C$, $\dag$ must also act as
the usual conjugation.  Furthermore, $\dag$ must also be compatible with
the other structures present i.e.
\begin{eqnarray}
(c_{X,Y})^\dag&=&(c_{X,Y})^{-1} \\
(\theta_X)^\dag&=&(\theta_X)^{-1}\\
(b_X)^\dag&=&d_Xc_{X,X^*}(\theta_X\otimes \one_{X^*})\\
(d_X)^\dag&=&(\one_{X^*}\otimes\theta_X^{-1})(c_{X^*,X})^{-1}b_X
\end{eqnarray}
\end{enumerate}
\begin{rmk} For any $f\in \Hom(X,Y)$ we define $f^*\in \Hom(Y^*,X^*)$
by:
$$f^*=
(d_Y\otimes \one_{X^*})(\one_{Y^*}\otimes f\otimes \one_{X^*})(\one_{Y^*}\otimes b_X).$$
\end{rmk}
\begin{rmk}
We will often consider categories satisfying some subset of these
axioms; for example a \emph{braided tensor category} satisfies
axioms 1-5.
\end{rmk}
\subsection{General Consequences}
The categorical axioms above supply us with several useful tools for studying
these categories.  The following results are found in the references mentioned
above or in \cite{OW} and \cite{TuW2}.

\subsubsection{Categorical Trace}
In any semisimple ribbon category one defines a categorical trace for
any morphism $f \in \End_\OO(X)$:
\begin{equation}
Tr_\OO(f)=d_X c_{X,X^*}((\theta_Xf)\otimes \one_{X^*})b_X: \1 \rightarrow \1.
\end{equation}
One defines the categorical dimension of an object $X$ by:
$$\dim_\OO(X):=Tr_\OO(\one_X).$$  It is often useful to normalize the trace
so that the trace of the identity morphism $\one_X$ has trace $1$ where $X$ is
any object.  This is achieved by setting $tr_\OO(f)=Tr_\OO(f)/\dim_\OO(X)$ for
any $f\in\End(X)$.

The expected properties of the trace go through and are by now well-known.
\begin{lemma}
(a)  $Tr_\OO(f\circ g)=Tr_\OO(g\circ f)$ when the composition and trace are
defined.\\
(b)  $Tr_\OO(f\otimes g)=Tr_\OO(f)Tr_\OO(g)$.\\
(c)  $\dim_\OO(X)\not=0$ if $X$ is simple.
\end{lemma}

A proof of the following important result can be found in \cite{OW}.
\begin{lemma}\label{parttrace}
 Let $\OO$ be a semisimple ribbon category, and $X$ and $Y$ be
 simple objects
in $\OO$, with $p\in \End(X\otimes X^*)$ the projection onto the
subobject of $X\otimes X^*$ isomorphic to $\1$, and $a\in
\End(Y\otimes X)$. Then
$$(\one_Y\otimes p)(a\otimes\one_{X^*})(\one_Y\otimes p)=
\frac{Tr_\OO(a)}{\dim_\OO(Y)\dim_\OO(X)}(\one_Y\otimes p)$$
\end{lemma}
The proof is an exercise in the so-called graphical calculus of
ribbon categories.  For an explicit formula for $p$ one may take
$1/\dim_\OO(X)b_Xb_X^\dag$ (which is defined regardless of the
existence of a conjugation in the category).

Lemma \ref{parttrace} has the following specialization known as
the \emph{Markov property} (see \cite{TuW2}):
\begin{lemma}\label{markov}
If $a\in\End(X^{\otimes n})$ and $m\in\End(X^{\otimes 2})$, then
$$tr((a\otimes \one_X)\circ (\one_X^{\otimes (n-1)}\otimes m))=tr(a)tr(m).$$
\end{lemma}
\subsubsection{Representations of $\C \B_n$}
The braiding axiom implies that the operators
$c_1:=c_{X,X}\otimes\one_X$ and $c_2:=\one_X\otimes c_{X,X}$ in
$\End_\OO(X^{\otimes 3})$ satisfy the braid relation
$c_1c_2c_1=c_2c_1c_2$, and hence we obtain representations of the
group algebra of the braid group $\C \B_n\rightarrow
\End_\OO(X^{\otimes n})$ by sending
$$\sigma_i \rightarrow
c_i:=\one_X^{\otimes (i-1)}\otimes c_{X,X}\otimes\one_X^{\otimes
(n-i-1)}$$ One may also define a representation of $\C \B_n$ on
the vector space $\End_\OO(X^{\otimes n})$ by composing with
$c_i$. Here $\sigma_i$ is the standard generator of $\B_n$ as
shown in Figure 1. \epsfverbosetrue \epsfxsize 2.2in
\begin{figure}[htb]
\centerline{\epsffile{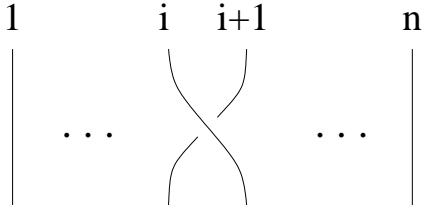}} \caption{The generator
$\sigma_i$}
\end{figure}

Tuba and Wenzl in \cite{TuW1} succeeded in classifying all
representations of $\B_3$ of dimension $\leq 5$ by the eigenvalues
of the image of $\sigma_1$ and the scalar by which the center of
$B_3$ acts.  This becomes quite useful as the structure of the
representation of $\C \B_n$ on $\End(X^{\otimes n})$ is already
essentially determined by considering $n=3$.

\subsubsection{Grothendieck Semiring}\label{grdef}
We also have the Grothendieck semiring $Gr(\OO)$ of a finite
semisimple monoidal category $\OO$.  Let $X_i$ $1\leq i\leq m$ be
a complete set of representatives of distinct isomorphism classes
of simple objects in $\OO$.  Define
$N^k_{ij}:=\dim\Hom(X_k,X_i\otimes X_j)$ so that $X_i\otimes
X_j\cong \bigoplus N^k_{ij} X_k$.  The braiding ensures that
$N_{ij}^k=N_{ji}^k$, and the Grothendieck semiring is the
commutative quotient ring:
\begin{equation}
Gr(\OO):=\Z[X_i]/\lan X_iX_j-\sum_kN^k_{ij}X_k: 1\leq i,j,k\leq
m\ra
\end{equation}
\begin{rmk}
The relations in $Gr(\OO)$ are often called \emph{fusion rules},
particularly in the physics literature.
\end{rmk}
The left action
of $Gr(\OO)$ on itself is given in the basis $\{X_i\}$ is
$X_i\rightarrow N_i$ where $(N_i)_{kj}=N_{ij}^k$. By commutativity
the simple representations of $Gr(\OO)$ are 1-dimensional, and we
can study the character theory.
\begin{definition} A \emph{character}
for $Gr(\OO)$ is any representation $f: Gr(\OO)\rightarrow
\C$ that satisfies
\begin{equation}\label{char}
f(X_i)f(X_j)=\sum_k N_{ij}^k f(X_k)
\end{equation}
\end{definition}
We have already seen one character of $Gr(\OO)$, namely the
function $\dim_\OO$.  Observe that for any character $f$ of
$Gr(\OO)$ the vector ${\bf f}:=\sum_i f(X_i)X_i$ is a simultaneous
eigenvector of the set of matrices $\M:=\{N_i\}$.  In particular
there can be at most $|\{X_i\}|$ inequivalent characters.

\subsection{Self-Dual Categories}\label{selfdual}
For convenience of notation, we make the following definition:
\begin{definition} A \emph{self-dual} category is one
in which every object is isomorphic to its dual object.
\end{definition}
All of the categories that we will consider in detail will be
self-dual. This gives $Gr(\OO)$ a much simpler structure:
 the $N_{ij}^k$ are totally symmetric in the $i,j$
and $k$. Lemma \ref{parttrace} has a stronger consequence in the
self-dual case (see \cite{TuW2}):
\begin{cor}\label{parttracecor}
Suppose $X^{\otimes 2}\cong \bigoplus_iX_i$ in a self-dual semisimple
ribbon category $\OO$, and we have a basis
of mutually annihilating idempotents $p_j\in\End(X^{\otimes 2})$ so that
$p_jX^{\otimes 2}\cong X_j$
and $X_1\cong \1$.
Then
$$(\one_X\otimes p_1)(p_j\otimes \one_X)(\one_X\otimes p_1)=
\frac{\dim_\OO(X_j)}{(\dim_\OO(X))^2}(\one_X\otimes p_1)$$
\end{cor}
\subsubsection{Unitary Categories}
In a semisimple Hermitian ribbon category, the form $\lan
f,g\ra:=Tr_\OO(f^\dag g)$ on $\Hom(X,Y)$ is Hermitian.  Since
$(c_{X,Y})^\dag=(c_{X,Y})^{-1}$ the form $\lan\ ,\ \ra$ on
$\End(X^{\otimes n})$ is preserved by the action of $\C \B_n$:
$$\lan c_if,c_ig\ra=Tr_\OO(f^\dag (c_i)^{-1}c_i g)=\lan f,g\ra.$$  So if
$\lan\ ,\ \ra$ is positive definite on $\End(X^{\otimes n})$ for
all $n$ then the representations of $\C \B_n$ is on the Hilbert
spaces $\End(X^{\otimes n})$ $n\geq 1$ are unitary.  For this
reason such categories are called \emph{unitary}.
\begin{definition}
A braided tensor category $\OO$ is called \emph{unitarizable} if
there is a Hermitian form on the morphism spaces with respect to
which the representations of $\C \B_n$ are unitary.
\end{definition}
If $\OO$ is a finite semisimple Hermitian ribbon
Ab-category then the positivity of the form $\lan\ ,\ \ra$ is
determined by positivity on the idempotents $p_i\in\End(X^{\otimes
n})$ where $p_iX^{\otimes n}\cong X_i$, since $\End(X^{\otimes
n})$ is a direct sum of full matrix algebras by semisimplicity.
Assume that $\OO$ is self-dual. Then we can choose the idempotents
so that $(p_i)^\dag=p_i$ (self-adjoint) and then we have that
$$\lan p_i,p_i\ra=Tr_\OO(p_i)=\dim_\OO(X_i)/(\dim_\OO(X))^n$$
by the lemmas above.  In particular, if $\dim_\OO(X_i)>0$ for all
simple objects $X_i$, then $\OO$ is unitary.
\begin{theorem}\label{nounitary}
Let $\OO$ be a semisimple self-dual Hermitian ribbon category.
Suppose that some simple object $X_i$ appears in $Y^{\otimes 2n}$
with $\dim_\OO(X_i)<0$. Then the category $\OO$ \emph{cannot} be
unitary for \emph{any} $\C \B_n$ invariant Hermitian form.
\end{theorem}
$Proof.$ Assume to the contrary that there is some $\C
\B_k$-invariant positive Hermitian form ( , ) on each
$\End(Y^{\otimes k})$.  Then choose positive projections
$p,p_i\in\End(Y^{\otimes 2n})$ with images isomorphic to $\1$ and
$X_i$ respectively.  By applying the corollary to Lemma
\ref{parttrace} in the self-dual case for $X=Y^{\otimes 2n}$ we
compute:
$$(\one_{Y}^{\otimes 2n}\otimes p)(p_i\otimes \one_Y^{\otimes 2n})
(\one_Y^{\otimes 2n}\otimes p)=
\frac{\dim_\OO(X_i)}{(\dim_\OO(Y))^{2n}}(\one_Y^{\otimes
2n}\otimes p).$$ But observe that the left-hand side of this
equality is a positive operator, whereas the right-hand side is a
negative operator as $\dim_\OO(X_i)<0$.\hfill $\Box$

\section{Ribbon Categories from Quantum Groups}
In this section we discuss the ribbon categories derived from the
quantum groups at roots of unity.  The construction of the
category is by now quite well-known and we will omit the details.
We content ourselves to refer the interested reader to: Jantzen's
book \cite{jantzen} for an introduction to quantum groups and to
illuminate the way through Lusztig's book \cite{lusztig} on the
same, the papers of Andersen and Paradowski \cite{andersen},
\cite{AndPar} for the categories of tilting modules and their
semisimple quotient, and chapters 9-11 of the book by Chari and
Pressley \cite{ChPr} for certain cases of the whole
construction.

\subsection{Notation and Outline}
We will need some notation in order to proceed.  Let $\g$ be a simple Lie
algebra of rank $k$.  We have:
\begin{itemize}
\item the root system $\Phi$ embedded in $\R^k$
\item the Cartan matrix $A=(a_{ij})$
\item a root basis $\Pi=\{\alpha_i\}_{i=1}^k$
\item positive roots $\Phi_+=\N\Pi\cap \Phi$
\item root lattice $Q=\Z\Pi$
\item a normalized bilinear form
$\lan\ ,\ \ra$ so that $2\lan \al_i,\al_j \ra/\lan \al_j,\al_j \ra=a_{ij}$
and $\lan \al,\al \ra=2$ for short roots.
\item coroot basis $\check{\Pi} =\{\check{\al_i} \}_{i=1}^k$,
$\check{\al_i}:=2\al_i/\lan \al_i,\al_i \ra$
\item coroot system $\check{\Phi}$
\item positive coroots $\check{\Phi}_+$
\item Weyl group $W$ generated by: $s_i(v):=v-\lan v,\check{\al_i}\ra\al_i$
\item fundamental weights $\La_i$ dual to the coroot basis via $\lan\ ,\ \ra$
\item weight lattice $P=\Z\{\La_i\}$
\item dominant Weyl chamber $C:=\R^+\{\La_i\}$ containing the dominant weights
$P_+:=\N\{\La_i\}$.
\end{itemize}

Drinfeld \cite{D} and Jimbo defined a quantum group as a
$q$-deformation $\mathcal{A}_q(\g)$ of the universal enveloping
algebra of $\g$ where the base field is $\Q(q)$ with $q$ an
indeterminate. The finite dimensional representations of
$\mathcal{A}_q(\g)$ are integral and highest weight and the
Grothendieck semiring of this representation category is
isomorphic to that of $\g$ itself. However, if we try to
specialize $q$ to a root of unity $\mathcal{A}_q(\g)$ is no longer
well-defined. Now let $q^2$ be a primitive $\ell$th root of unity,
that is, $q=e^{z\pi i/\ell}$ with $\gcd(z,\ell)=1$.  Lusztig's
``modified form" of $\mathcal{A}_q(\g)$ denoted
 $U_q\g$ is well-defined for any complex
$q\not\in \{-1,0,1\}$.  That $U_q\g$ is a ribbon Hopf algebra
follows from the work of Drinfeld, Lusztig and others, see
\cite{BK} for details. For each $\la\in P_+$ one constructs a
\emph{Weyl module} $V_\la$ of $U_q(\g)$ by restricting the
corresponding highest weight $\mathcal{A}_q(\g)$-module to
$U_q(\g)$ and specializing the parameter $q$ to the chosen root of
unity. The generators of $U_q\g$ act on Weyl modules by matrices
with entries in $\Z[q^{\pm 1}]$.  The Weyl modules are not all
irreducible or even indecomposable. To remedy this, Andersen
\cite{andersen} defined a category $\T$ of \emph{tilting modules}
that have the following key properties:
\begin{enumerate}
\item[1.] For each $\la\in P_+$ there is a unique indecomposable tilting
module $T_\la$.
\item[2.] The set $\I=\{T\in \T: \dim_\T(T)=0\}$ is a tensor
ideal.
\item[3.] There are finitely many indecomposable $T_\la\not\in
\I$.  These are irreducible and isomorphic to the corresponding
Weyl module.
\item[4.] The category $\F=\T/\I$ defined by taking the
quotient of the morphisms in $\T$ by the radical of $Tr_\T$ is a
semisimple ribbon category.
\end{enumerate}

Tilting modules can be realized as
direct sums of submodules of tensor powers of the
\emph{fundamental module(s)}.  A fundamental module is one that
generates the category $Rep(U_q\g)$ generically, that is every
irreducible module appears in some tensor power.

\subsection{The Category $\F$}
We can describe the category $\F$ as
follows.  Let $d$ be the ratio of the squared length of a long
root in $\g$ to the squared length of a short root.  If $q^2$ is a
primitive $\ell$th root of unity with $\ell$ divisible by $d$ then let
$\theta$ be the highest root of $\Phi$, if $\ell$ is coprime to
$d$ then let $\theta$ be the highest short root.  Then the simple
objects of $\F$ are isomorphic to Weyl modules $V_\la$ with
$\la\in C_\ell:=\{\mu\in P_+: \lan \mu+\rho,\check{\theta}\ra < \ell\}$, here
$\rho$ is half the sum of the positive roots $\al\in\Phi_+$.
In fact the indecomposable tilting modules $T_\mu$ that are isomorphic to
irreducible Weyl modules
are labeled by
$$\mu\in\overline{C_\ell}:=\{\mu\in P_+: \lan \mu+\rho,\check{\theta}\ra \leq \ell\}.$$
To avoid degeneracies we always assume that the rank $k$ and
$\ell$ are such that $\rho+\La_1\in C_\ell$, where $\La_1$ is the
dominant weight of the defining representation of $\g$.  By taking
the convex hull of the set $C_\ell$ we obtain the
\emph{fundamental Weyl alcove} denoted by $D$.

\subsubsection{Affine Weyl Group}
The dominant Weyl chamber $C$ is described as the fundamental
domain of the Weyl group $W$ containing $\rho$, the fundamental
Weyl alcove $D$ can be similarly described:
\begin{definition}
Denote the affine reflection in $\R^k$ through the hyperplane
$\{x\in \R^k: \lan x, \check{\theta}\ra=\ell\}$ by $t_\ell$.  If
we adjoin
 $t_\ell$ to the Weyl group $W$ we get the \emph{affine}
Weyl group $W_\ell$.  Explicitly $t_\ell(\la)=\la+(\ell-\lan
\la,\check{\theta} \ra)\theta.$
\end{definition}
We must define a slightly different action of $W_\ell$ on $P$ than
the usual one inherited from Euclidean $\R^k$.  For $w\in W_\ell$
and $s\in \R^k$ define the ``dot action" $w\cdot
x:=w(x+\rho)-\rho$. Then $D$ is the fundamental domain of the dot
action of $W_\ell$ on $\R^k$ and of course $C_\ell=D\cap P_+$. The
elements of $W$ have a natural signature $\sgn$ depending on the
number of simple reflections $s_i$ in any decomposition.  If we
assign $\sgn(t_\ell)=-1$ then this extends the signature function
to $W_\ell$.

We now proceed to
describe the categorical structure on $\F$.
\subsubsection{Monoidal Structure}
$\F$ inherits a monoidal structure from the
comultiplication and counit in the Hopf algebra $U_q\g$.
\subsubsection{Duality} The dual module of a simple Weyl module
$V_\la$ is the ordinary vector space dual with the action of
$U_q\g$ defined via the antipode.  $V_\la^*$ is also a Weyl module
with highest weight equal to $-w_0(\la)$ where $w_0$ is the
longest element in the Weyl group with respect to Bruhat order.
One checks that $-w_0(\la)\in C_\ell$.  The rigidity morphisms are
defined
$$b_{V}: 1\rightarrow \sum_i v_i\otimes v^i$$
and
$$d_{V}: f\otimes v\rightarrow f(v)$$
where $v_i$ is a basis of $V$ and $v^i$ is the dual basis (of
$V^*$).
\subsubsection{Braiding}\label{rmat}
Lusztig \cite{lusztig} showed that the universal $R$-matrix in
$U_q\g$ specializes to the root of unity case.  Composing with the
flip operator $\sigma$ we get well-defined operators
$\check{R}_{V,W}$ for any objects $V,W\in\F$.  These do satisfy
the braiding identities. We have the very useful (see \cite{D}):
\begin{prop}[Drinfeld]\label{reigs}
If $V_\la$ and $V_\mu$ are simple Weyl
modules such that $V_\nu$ appears in $V_\la\otimes V_\mu$ then one
has:
$$\check{R}_{\mu,\la}\check{R}_{\la,\mu}\mid_{V_\nu}=q^{c_\nu-c_\la-c_\mu}\one_{V_\nu}$$
where $c_\gamma:=\lan\gamma+2\rho,\gamma\ra$.
\end{prop}
To conform with our original notation we will denote the morphisms
$\check{R}_{V,W}$ by $c_{V,W}$.
\subsubsection{Twist}
It also follows from the work of Drinfeld that there is a
\emph{universal Casimir} operator in $U_q\g$ that provides $\F$
with a twist.  For a simple object $V_\la$ the twist $\theta_\la$
acts by the constant $q^{c_\la}$ where $c_\la$ is as above.
\subsubsection{Ab-structure}
The spaces $\Hom(V,W)$ are quotients of the vector spaces of
intertwining operators in the category $\T$, so
they are themselves $\C$-vector spaces.
\subsubsection{Finite Semisimplicity}
Andersen's (\cite{andersen}) main result shows that $\F$ is
semisimple, as we have taken the quotient by the radical part of
the category, and all other necessary properties are inherited
from the category $\T$.  Only finitely many isomorphism classes of simple objects
$V_\la$ survive in the quotient.
\subsubsection{Hermitian Form}
Kirillov Jr. \cite{kirillov} succeeded in defining a conjugation on the
category $\F$.  In this paper we are only concerned with the
existence of one, so we will not go into details.
\subsubsection{Categorical Trace}\label{qdim}
With all of the above structure, $\F$ is a ribbon category and
hence has a trace.  We can compute the value of $\dim_\T$
explicitly on the objects $T_\mu\cong V_\mu$,
$\mu\in\overline{C_\ell}$:
$$\dim_\T(V_\mu)=\prod_{\al\in\Phi_+}\frac{[\lan\mu+\rho,\al\ra]}{[\lan\rho,\al\ra]}$$
where $[n]:=\frac{q^n-q^{-n}}{q-q^{-1}}$.  This follows from the
proof of the Weyl dimension formula in the classical theory.
Since $\theta\in\Phi_+$ one sees that $\dim_\T(V_\mu)=0$ for
$\mu\in(\overline{C_\ell}\setminus C_\ell)$.  By construction $\dim_\T$
vanishes on the ideal $\I$ so the categorical dimension
$\dim_\F$ coincides with $\dim_\T$ on the quotient.
\subsubsection{Grothendieck Semiring}
The Grothendieck semiring $Gr(\F)$ is a quotient of $Gr(Rep(\g))$.
The structure constants of $Gr(\F)$ are
$W_\ell$-antisymmetrizations of
 those of $Gr(Rep(U_q\g))$ for
$q$ generic (which are the same as those of $Gr(Rep(\g))$).
\begin{rmk} The proposition that follows was
proved for weights in the \emph{root} lattice by Andersen and
Paradowski (\cite{AndPar}, Prop. 3.20), as the quantum group
studied there is constructed from the \emph{adjoint} root datum
whereas we want to use the \emph{simply connected} root datum (see
\cite{lusztig}, Chapter 2).  So in particular one must justify the
extension of this result to those weights not in the root lattice, that is, the
half-integer weights.  However, the argument presented in \cite{AndPar} relies only upon results in \cite{andersen} (which are valid for $\ell$ coprime to the nonzero entries of the Cartan matrix, in particular for Lie type B quantum groups with $\ell$ odd: see Section 1 of \cite{andersen}) and therefore carries over
word-for-word to the case at hand. In fact, the only results cited by Andersen and Paradowski for which they do not give an explicit reference in \cite{andersen} are the linkage principle and their ``quantum version of Proposition 2.5''.  These are found in \cite{andersen} statement (1.2) and Theorem 2.5 respectively.  With this justification we attribute the proposition below to Andersen and Paradowski.\end{rmk}

If $m_{\la\mu}^{\nu}=\dim\Hom_{U_q\g}(V_\nu,V_\la\otimes V_\mu)$
for $\la,\mu,\nu\in P_+$ (that is, $m_{\la\mu}^{\nu}$ are the
classical weight multiplicities), then we have (see
\cite{AndPar}, Prop. 3.20, see also \cite{S2} for the general case):
\begin{prop}[Andersen-Paradowski]\label{antisym}
For simple objects $V_\la, V_\mu$ in the category $\F$,
$$N_{\la\mu}^{\nu}=\sum_{w\in
W_\ell:\quad w\cdot\nu\in P_+}\ve(w)m_{\la\mu}^{w\cdot\nu}$$ where
$N_{\la\mu}^{\nu}:=\dim\Hom_\F(V_\nu,V_\la\otimes V_\mu)$.
\end{prop}

Observe that if $\lan\nu+\rho,\check{\theta}\ra=\ell$
(i.e. $\nu\in\overline{C_\ell}\setminus C_\ell$ and
$\dim_\T(V_\nu)=0$) then
$$t_\ell\cdot\nu=t_\ell(\nu+\rho)-\rho=
(\nu+\rho)+(\ell-\lan (\nu+\rho),\check{\theta}
\ra)\theta-\rho=\nu$$ so the antisymmetrization above gives
$N_{\la\mu}^{\nu}=0$ as expected.

\section{Type $B$ at Odd $\ell$ Categories}\label{typeB}
Observe that the construction of the categories above depend on
two choices: a Lie algebra $\g$ and a root of unity $q^2$.   We
now specialize to the categories we will study in detail: that is,
the Lie algebra $\g\cong\mathfrak{so}_{2k+1}$, and $q^2$ a
primitive $\ell$th root of unity, $\ell$ odd.  For a fixed $\ell$
and $k$, we denote by $\mF$ the family of ribbon categories
constructed as above from $\mathfrak{so}_{2k+1}$ with $q^2$ any
primitive $\ell$th root of unity.  A fixed member of this family
will be denoted by $\F$.

\subsection{Type $B$ Data}
Let
$\{\ve_i\}$ be the standard basis for $\R^k$.  We fix a root basis
$$\Pi=\{\al_i\}_1^k =\{\ve_1-\ve_{2},\ve_2-\ve_3,\ldots,\ve_{k-1}-\ve_{k},\ve_k\}$$
so the root lattice $Q=\rm{span}_\Z\{\al_i\}_1^k$ is just
$\Z^k$.  We also record that the set of positive roots is
$$\Phi_+=\{\ve_s\pm\ve_t,\ve_u: s<t\}.$$
The form $\lan\ ,\ \ra$ is twice the usual
dot product on $\R^k$ so that the square length of long roots is
$4$, and $2$ for short roots. Thus the coroot basis
$\check{\Pi}=\{\check{\al}\}$ has
$$\check{\al}_i=\begin{cases}\frac{1}{2}(\ve_i-\ve_{i+1}) & i=1,\ldots,k-1\\
\ve_k & i=k
\end{cases}$$
Note that classically, the coroots for type B are the roots of type
C, but here we must take care as the
normalization of the form is not the classical
one.  We will see where this leads to subtleties later.  The Weyl
group $W$ is the semi-direct product of $S_k$ and $(\Z_2)^k$ and
acts on $\R^k$ via permutations and sign changes.

For our choice of a root basis we
have the following {\it fundamental weights}:
$$\La_i=\begin{cases}\sum_{1\leq j \leq i-1}\ve_j & i \leq k-1\\
\frac{1}{2}\sum_{1\leq j \leq k}\ve_i & i=k \end{cases}$$
and the
dominant weights $P_+= {\rm span}_{\N}\{\La_j\}_1^k$.  The weight lattice $P={\rm
span}_\Z\{\La_j\}_1^k$ is then seen to be $\Z^k \bigcup
(\La_k+\Z^k)$. For convenience of notation we introduce the
function on $P$:
$$p(\la)=\begin{cases} 1 & \text{if $\la \in \Z^k$}\\
-1 & \text{if $\la \in (\La_k+\Z^k)$}
\end{cases}$$  We refer to a weight $\la$ as integral, resp. half-integral,
if $p(\la)=1$, resp. $p(\la)=-1.$
The weights are usually represented as $k$-tuples,
e.g. $\La_k=(1/2,\ldots,1/2)$.

For type $B$ we have that
$w_0=-1$, that is, the element of the Weyl group that changes the sign of each
coordinate.  Since the weight of $V_\la^*$ is $-w_0(\la)=\la$, all modules are
self-dual in the present case.

\subsection{Classical Representation Theory,
Abridged}\label{tensorreps} As we noted above, for generic
parameters $q$, we have $Gr(Rep(U_q\g))\cong Gr(Rep(\g))$, and
$Gr(\F)$ is a quotient of these rings, so in this subsection we
will summarize the necessary facts from the representation theory
of the algebra $U\mathfrak{so}_{2k+1}$. This material can be found
in any introductory text on Lie groups, such as \cite{bigred} or
\cite{humphreys}, and goes back at least to Weyl \cite{Wy}.

The irreducible finite-dimensional integral highest weight modules of
$\mathfrak{so}_{2k+1}$ are in one-to-one correspondence with the
elements of $P_+$.  Each irreducible integral highest weight module
 $V_\la$ has a multiset
of weights $P(\la)$ which correspond to the weight-space
decomposition of $V_\la$
with respect to the action of the Cartan subalgebra.  The multiset
$P(\la)$ lies in the ball of radius $|\la|$ (ordinary
Euclidean distance) centered at the origin, and the weights
in the $W$-orbit of $\la$ appear with multiplicity one.
The other weights are of the form $\la-\al$ for some
$\al \in \Q$.
To decompose the tensor product of two irreducible
modules $V_\la$ and $V_\mu$ one looks at the
intersection
$\{\nu=\mu+\kappa: \kappa \in P(\la)\} \bigcap P_+$ which contains
the dominant weights of the irreducible submodules
$$P_+(V_\la \otimes V_\mu)=\{\nu\in P_+:V_\nu\subset V_\la\otimes V_\mu\}.$$
We do not formulate the
precise algorithm to determine which $V_\nu$ do occur nor the multiplicities, but
we can say that the irreducible module
$V_{\mu+w(\la)}$ appears with multiplicity one, where $w$ is any element
in the Weyl group such that $w(\la)+\mu \in P_+$.  (This follows from
the outer multiplicity formula, see e.g. \cite{bigred} Corollary 7.1.6).  Moreover,
$P_+(V_\la \otimes V_\mu)$ is contained in the ball of radius $|\la|$ centered
at $\mu$, and $p(\nu)=p(\la)p(\mu)$ for any $\nu \in P_+(V_\la\otimes V_\mu)$.
In other words, all weights of
simple submodules of $V_\la\otimes V_\mu$ are integral if $\la$ and
$\mu$ are both integral or half-integral, and half-integral otherwise.

\subsection{Structure Constants of $Gr(\F)$}
Recall the left regular representation of $Gr(\F)$ from
\ref{grdef} and denote the images of the generators by $N_\la$,
$\la\in C_\ell$.
In general it is not easy to compute the entries $N_{\la\mu}^{\nu}$ of
the matrices $N_\la$ as
it is already difficult to compute the classical
multiplicities $m_{\la\mu}^\nu$;
however, for our analysis we only require
two explicit decomposition rules--both of which
were already known to
Brauer in the 1940s.  We begin with the
decomposition rules for tensoring with the generating module $V_{\La_k}$.
\begin{example}\label{minrules}
We have that $V_{\La_k}$ is a minuscule representation (all weights are conjugate
under the Weyl group) the simple decomposition as a
$\mathfrak{so}_{2k+1}$-module is:
$$V_{\La_k}\otimes V_\la\cong\bigoplus_{W^\la_k}
V_{\la+w(\La_k)}$$ where $W^\la_k=\{ w\in W:
\la+w(\La_k)\in P_+\}$
Note that $W(\La_k)=\{\frac{1}{2}(\pm 1,\ldots,\pm 1)\}$, so
all $\la+w(\La_k)$ are in $\overline{C_\ell}$, so the $W_\ell$-antisymmetrization
has the effect of discarding the $V_{\la+w(\La_k)}
\in \overline{C_\ell}\setminus C_\ell$ and leaving all other objects alone.
That is, for $\la,\nu \in C_\ell$
\begin{equation}\label{rule1}
N_{\La_k\la}^\nu=\begin{cases}1 & \text{if $\nu=\la+w(\La_k)$ some $w \in W$}\\
0 & \text{otherwise}
\end{cases}
\end{equation}

\begin{lemma}\label{Fgen}
$V_{\La_k}$ generates $\F$.
\end{lemma}
$Proof.$
We will show that there exists an odd integer $s$
such that \emph{every} simple object in $\F$ appears in $V^{\otimes s}_{\La_k}$ or
$V_{\La_k}^{\otimes s+1}$.  Every weight $\la\in C_\ell$ can be expressed as
a sum of weights in $W(\La_k)$, so every $V_\la$ appears
in some tensor power of $V_{\La_k}$ by an induction using the multiplicity formula
above.  Furthermore,
the trivial representation $\1$ appears
in $V_{\La_k}^{\otimes 2}$ so once $V_\la$ appears in an odd (resp. even)
tensor power of $V_{\La_k}$ it will appear in \emph{every} odd (resp. even)
tensor power thereafter.\hfill $\Box.$
\end{example}

The vector (or defining) representation of $\mathfrak{so}_{2k+1}$
has
highest weight $\La_1=\ve_1$.  We will only need to know the
decomposition for tensoring $V_{\La_1}$ with simple objects whose highest
weights have integer entries:
\begin{example}\label{rule2}
The weights of $V_{\La_1}$ are the zero weight together with
$W(\La_1)=\{\pm\ve_i: 1 \leq i \leq k\}$.  The decomposition
algorithm as a $\mathfrak{so}_{2k+1}$-module is (for integral weights $\mu$):
$$V_{\La_1}\otimes V_\mu\cong\delta(\mu) V_\mu+\bigoplus_{W_1}V_{\mu+w(\La_1)}$$
where $W_1=\{w \in W: w(\La_1) \in P_+\}$ and $\delta(\mu)=1$ if
$\lan \mu,\ve_k \ra >0$ and zero otherwise.  Since the dominant weights
in $\overline{C_\ell}\setminus C_\ell$ all have integer entries and
$\mu$ is distance at least 1 from the hyperplane spanned by
$\overline{C_\ell}\setminus C_\ell$ we conclude that
$P_+\cap(\mu+W(\La_1)) \subset \overline{C_\ell}$.  Hence the
$W_\ell$-antisymmetrization has the effect of
discarding those $V_{\mu+w(\La_1)}$ with
$\mu+\La_1\in\overline{C_\ell}\setminus C_\ell$.
So for $\mu,\nu \in C_\ell\cap\Z^k$ we compute:
\begin{equation}
N_{\La_1\mu}^\nu=
\begin{cases}1 & \text{if $\nu=\mu \pm \ve_i$ for some $1\leq i\leq k$}\\
1 & \text{if $\mu=\nu$ and $\lan \mu,\ve_k\ra > 0$}\\
0 & \text{otherwise}
\end{cases}
\end{equation}
As in Example \ref{minrules}, we can use this computation to
conclude that $V_{\La_1}$ generates the subcategory of objects
labeled by integer weights.  It is slightly trickier to show that,
in fact, every object labeled by an integer weight appears in both
an even \emph{and} an odd power of $V_{\La_1}$.  The trick is to
find a $\mu$ with $|\mu|=s$ odd and $\lan \mu,\ve_k\ra > 0$.  Then
$V_\mu$ first appears in $V_{\La_1}^s$ (that is, $s$ is minimal
with this property).  But then $V_\mu$ appears in
$V_{\La_1}^{s+1}$ by the rule above.  Since $s+1$ is even, by
applying our rule again and using the fact that $N_{\La_1
\mu}^{\nu}$ are completely symmetric we see that $V_{\La_1}$
appears in $V_{\La_1}^{(s-1)(s+1)}$. Thus \emph{every} object
labeled by an integer weight appears in an even tensor power of
$V_{\La_1}$.  By considering cases and applying this argument
again we conclude that the same is true for odd tensor powers of
$V_{\La_1}$.
\end{example}

\subsection{Character Analysis}
Eventually we want to describe all (irreducible) characters of the ring $Gr(\F)$.
Our main source of characters are the \emph{$q$-characters} of
$Gr(Rep(\mathfrak{so}_{2k+1}))$, which are nothing more than
ordinary
characters of the ring $Gr(Rep(U_q\mathfrak{so}_{2k+1}))$ for $q$ generic.
To each $\nu\in Q$ (the root lattice) there is a certain ``diagonalizable'' element in
the quantum group $U_q(\mathfrak{so}_{2k+1})$ denoted by $H_\nu$ on which the
irreducible characters $\chi_\la$ are defined for any $\la\in P_+$.  This is completely
analogous to the classical situation where the characters act on a maximal abelian
subalgebra.
\begin{equation*}
\chi_\la(H_\nu)= \frac{1}{\delta_B(H_\nu)}
\sum_{w \in W}\ve(w)q^{\lan w(\la+\rho),\nu \ra}
\end{equation*}
where $$\delta_B(H_\nu)=\sum_{w \in W}\ve(w)q^{\lan w(\rho),\nu
\ra}$$
is the \emph{Weyl denominator}.  Recall that $[n](q-q^{-1})=q^n-q^{-n}$.  An
important computation due to Weyl \cite{Wy} gives
us the product form
$$\delta_B(H_\nu)=\prod_{\al \in
\Phi_+}[\frac{1}{2}\lan \al,\nu \ra]$$
(see \cite{bigred} Chapter 7 for
a more modern treatment).  The $\frac{1}{2}$ appears
here because we have normalized the form $\lan\  ,\  \ra$ to be
twice the form used in the classical theory. (Note that
$\frac{1}{2}\lan\al,\nu \ra$ is a integer since both $\al$
and $\nu$ have integer entries.)

For any fixed $\nu\in Q$ the characters $\chi_\la$ satisfy:
\begin{enumerate}
\item[1.] $\chi_\zero(H_\nu)=1$
\item[2.]
$\chi_\la(H_\nu)\chi_\mu(H_\nu)=\sum_\kappa m_{\la\mu}^\kappa\chi_\kappa(H_\nu)$
where
$$m_{\la\mu}^\kappa=\dim\Hom_{U_q\mathfrak{so}_{2k+1}}(V_\kappa,V_\la\otimes V_\mu)$$
\end{enumerate}
The first property is clear, while the second is a fundamental result in
classical representation theory.

Now suppose $q^2$ is a primitive $\ell$th root of unity.  Notice
 that if $\nu=2\rho$ and
$\la\in \overline{C_\ell}$
Weyl's formula give us:
$$\chi_\la(H_{2\rho})=\dim_\T(V_\la).$$
This motivates the following notation:
\begin{definition}
Let $\mu \in P_+\setminus \Z^k$
so that $\mu+\rho \in P_+\cap Q$ (i.e. $p(\mu+\rho)=1$).
Then for all $\la \in
P_+$ we define
$$\dim_\F^{\mu}(V_\la):=\chi_\la(H_{\mu+\rho}).$$
\end{definition}
The following technical lemma gives the precise criterion for a
character of $Gr(Rep(U_q(\mathfrak{so}_{2k+1})))$ to specialize to
a character of $Gr(\F)$:
\begin{lemma}
The specialization of a character $\chi_\kappa(H_\nu)$ to
$Gr(\F)$ gives a character of $Gr(\F)$ if and only if:
\begin{enumerate}
\item[3.] $\chi_\kappa(H_\nu)=\ve(w)\chi_{w\cdot\kappa}(H_{\nu})$
for all $\kappa \in C_\ell$, all $w
\in W_\ell$ such that $w\cdot\kappa\in P_+$ and $q^2$ an $\ell$th root of unity,
$\ell$ odd.
\end{enumerate}
\end{lemma}
$Proof.$
Setting
$W_\kappa=\{w \in W_\ell: w\cdot\kappa \in P_+\}$ for $\kappa \in \overline{C_\ell}$,
the second property of characters $\chi_\la$ becomes:
\begin{align*}
\chi_\la(H_\nu)\chi_\mu(H_\nu)=\sum_im_{\la\mu}^{\mu_i}\chi_{\mu_i}(H_\nu)&=\\
\sum_{\kappa \in \overline{C_\ell}} \left(\sum_{w \in
W_{\kappa}}\ve(w)m_{\la\mu}^{w\cdot\kappa}\right)\chi_{\kappa}(H_\nu)&= \sum_{\kappa \in
C_\ell}N_{\la\mu}^{\kappa}\chi_{\kappa}(H_\nu)
\end{align*}
since
to every $\mu_i \in P_+$ there is a unique
$\kappa \in \overline{C_\ell}$ so that $w\cdot\kappa=\mu_i$ for some
$w \in W_\ell$ and $N_{\la\mu}^{\kappa}=0$ if
$\kappa \in \overline{C_\ell}\setminus C_\ell$.
\hfill $\Box$\\

To prove Property 3 in the above lemma
we need only verify it for simple reflections $s_i$,
$t_\ell$ since they generate
$W_\ell$.  Moreover, we need only consider the numerator of
$\chi_\kappa(H_\nu)$ as the denominator $\delta_B(H_\nu)$ does not
depend on $\kappa$. So the veracity of Property 3
 will follow from the following lemma:
\begin{lemma}
$\sum_{w \in W}\ve(w)q^{\lan w(r\cdot\kappa+\rho),\nu \ra}=
\ve(r)\sum_{w \in W}\ve(w)q^{\lan
w(\kappa+\rho),\nu \ra}$ for $r$ a simple reflection and $\nu\in Q$.
\end{lemma}
$Proof.$
Define $w^{\prime} \in W$ by $w^{\prime}(\la)=\la-\lan \la,\ve_1\ra\ve_1$
and observe that $\ve(w^{\prime})=-1$ as $w^{\prime}$ just changes the
sign of the first coordinate of $\la$.
We compute:
\begin{align*}
&\lan w(t_\ell\cdot\kappa+\rho),\nu \ra=
\lan t_\ell(\kappa+\rho)-\rho+\rho,w^{-1}(\nu) \ra=\\
&\lan (\kappa+\rho)-\lan\kappa+\rho,\ve_1\ra\ve_1+\ell\ve_1,w^{-1}(\nu) \ra=
\lan ww^\prime(\kappa+\rho),\nu\ra+\ell\lan\ve_1,\nu\ra
\end{align*}
Since $\ell\lan\ve_1,\nu\ra$ is an even multiple of $\ell$ and $\ve(t_\ell)=-1$,
we have:
\begin{equation*}
\sum_{w \in W}\ve(w) q^{\lan w(t_\ell\cdot\kappa+\rho),\nu \ra}=
\ve(t_\ell)\sum_{w \in W}\ve(w)q^{\lan
w(\kappa+\rho),\nu \ra}
\end{equation*}
after reindexing the sum.
The computation for $s_i$ is slightly
less complicated, and just follows from the fact that
$\chi_\kappa(H_\nu)$ is an antisymmetrization with respect to the
Weyl group of the characters of the finite abelian group $\ell P/Q$.  It can
also be computed directly as for $t_\ell$.
Thus we have proved the lemma. \hfill $\Box$\\
Thus the specialization to roots of unity and restriction to $C_\ell$ of the
characters $\chi_\kappa(H_\nu)$ are indeed characters of the ring $Gr(\F)$.

Next we prove the following crucial:
\begin{lemma}\label{pos}
$\dim^{\La_k}_\F(V_{\la})$ is positive for all $\la \in C_\ell$
for $q=e^{\pi i/\ell}$.
\end{lemma}
$Proof.$
First we consider the numerator
$$\sum_{w \in W}\ve(w)q^{\lan w(\la+\rho),\La_k+\rho\ra}$$ of
$\dim^{\La_k}_\F(V_{\la}).$  Observe that the positive coroots
$\check{\alpha} \in \check{\Phi}_+$ are $\frac{1}{2}$
the positive roots $\Phi_+^C$ of type $C$ (corresponding to
$\mathfrak{sp}_{2k}$). In the classical theory we would get
exactly the positive roots of type $C$, but we are using twice the
classical form. Furthermore $\La_k+\rho=\rho^\prime$ is one-half
the sum of the positive roots of type $C$ and is thus the sum of
the positive coroots as we have defined them. Moreover, the Weyl
group $W$ is the same for these two algebras.  Let ( , ) be the
usual inner product on Euclidean space, so that $2(a,b)=\lan a,b
\ra$.  We have that
\begin{align*}
\sum_{w \in W}\ve(w)q^{\lan w(\la+\rho),\La_k+\rho\ra}&=
\sum_{w \in W}\ve(w)q^{\lan \la+\rho, w(\rho^\prime)\ra}
=\sum_{w \in W}\ve(w)q^{(2(\la+\rho),w(\rho^\prime))} \\
=\prod_{\beta \in \Phi^C_+}[(\la+\rho,\beta)]
&=\prod_{\check{\alpha} \in \check{\Phi}_+}[2(\la+\rho,\check{\alpha})]
=\prod_{\check{\alpha} \in \check{\Phi}_+}[\lan \la+\rho,\check{\alpha}\ra]
\end{align*}
by the observations above and the classical Weyl denominator
factorization for type $C$.  The same computation for $\la=\zero$
shows that the denominator of $\dim^{\La_k}_\F(V_{\la})$ also
factors nicely so that when we evaluate at $q=e^{\pi i/\ell}$ we
get:
\begin{align*}
\dim^{\La_k}_\F(V_{\la})&=\prod_{\check{\alpha} \in
\check{\Phi}_+} \frac{[\lan
\la+\rho,\check{\alpha}\ra]}{[\lan\rho,\check{\alpha}\ra]}
=\prod_{\check{\alpha} \in \check{\Phi}_+} \frac{q^{\lan
\la+\rho,\check{\alpha}\ra}-q^{-\lan \la+\rho,\check{\alpha}\ra}}
{q^{\lan \rho,\check{\alpha}\ra}-q^{-\lan \rho,\check{\alpha}\ra}} \\
&=\prod_{\check{\alpha} \in \check{\Phi}_+}
\frac{\sin(\lan \la+\rho,\check{\alpha}\ra\pi i/\ell)}
{\sin(\lan\rho,\check{\alpha}\ra\pi i/\ell)}.
\end{align*}
Now we see
that when $\la \in C_\ell$, $\lan \la+\rho,\check{\alpha}\ra<\ell$ for all
$\check{\alpha} \in \check{\Phi}_+$ so that each factor in
the above product is positive.\hfill $\Box$

We end this section with an important uniqueness theorem which relies on
the classical theorem of Perron and Frobenius found in \cite{Ga}.  Recall
that a positive matrix is a matrix whose entries are all strictly positive.
\begin{prop}[Perron-Frobenius]
A positive matrix $A$ always has a positive real eigenvalue of multiplicity
one whose modulus exceeds the moduli of all other eigenvalues.  Furthermore
the corresponding eigenvector may be chosen to have only positive real
entries and is the unique
eigenvector with that property.
\end{prop}
We now proceed to prove:
\begin{theorem}\label{Unique}
Evaluating $\dim^{\La_k}_\F(V_{\la})$ at $e^{\pi i/\ell}$ gives
the only character of $Gr(\F)$ that is positive for all $\la
\in C_\ell$.
\end{theorem}
$Proof.$ We observed in \ref{selfdual} that if $f:
C_\ell\rightarrow \C$ is a character of $Gr(\F)$ then the vector
${\bf f}=(f(\la))_{\la \in C_\ell}$ must be a simultaneous
eigenvector of the set $\mathcal{M}:=\{N_\la\}$, $\la\in C_\ell$.
In fact, using the definition of $N_\la$ one computes that
$N_\la({\bf f})=f(\la){\bf f}$. So if we can show that $N_{\La_k}
\in \mathcal{M}$ has only one positive eigenvector we will have
proved the theorem.  In the proof of Lemma \ref{Fgen} we saw that
for some odd integer $s$, the matrix $N_{\La_k}^s+N_{\La_k}^{s+1}$
has all positive entries. So one may apply the  Perron-Frobenius
Theorem to the matrix $N_{\La_k}^s+N_{\La_k}^{s+1}$ to see that it
has a {\it unique} positive eigenvector.  But $N_{\La_k}$ is a
(symmetric) diagonalizable matrix, so it has the same eigenvectors
as $N_{\La_k}^s+N_{\La_k}^{s+1}$. Since $\dim^{\La_k}_\F(V_{\la})$
at $e^{\pi i/\ell}$ was shown to be positive in Lemma \ref{pos},
we are done. \hfill $\Box$

\subsection{The Involution}\label{invol}
Next we define an involution $\phi$
of $C_\ell$ that will be central to the analysis of the characters
of $Gr(\F)$.
Let $\ga \in C_\ell$ be such that $|\ga|$ is maximal, explicitly,
$\ga=(\tfrac{\ell-2k}{2},\ldots,\tfrac{\ell-2k}{2}).$  Further denote by $w_1$ the
element of the Weyl group $W$ such that
$w_1(\mu_1,\ldots,\mu_k)=(\mu_k,\ldots,\mu_1).$
Define $\phi(\la):=\ga-w_1(\la)$.  It is clear that $\phi$ is a bijective map
from $C_\ell$ to itself and
that $\phi^2(\la)=\la$, and that $\phi \not\in W_\ell$ as no $\la \in P_+$ is
fixed by $\phi$.  The following lemma describes the key property of
$\phi$.
\begin{lemma}\label{phi}
For $q^2$ a primitive $\ell$th root of unity the involution $\phi$
preserves $|\dim_\F^{\mu}|$ (for $\mu \in P_+\setminus \Z^k$),
that is,
\begin{equation}
\dim_\F^{\mu}(V_\la)=\pm\dim_\F^{\mu}(V_{\phi(\la)})
\end{equation}
In particular (by setting $\mu=\rho$) this holds for the
categorical dimension $\dim_\F$ of $\F$.
\end{lemma}
$Proof.$
 Fix $\mu \in P_+\setminus \Z^k$
and a choice of a primitive $\ell$th root of unity $q^2$ (so
$q^\ell=\pm 1$). First consider $\sum_{w \in W}\ve(w)q^{\lan
\la+\rho,w(\mu+\rho) \ra}$ the numerator of
$\dim_\F^{\mu}(V_\la)$.  We compute
\[
\begin{split}
\lan \phi(\la)+\rho,w(\mu+\rho)\ra &=\lan
\ga-w_1(\la)+\rho,w(\mu+\rho)\ra \\
&=\lan w_1(\ga-\la+\rho+w_1(\rho)-\rho),w(\mu+\rho)\ra \\
&=\lan \ga+\rho+w_1(\rho),w_1 w(\mu+\rho)\ra +
\lan \la+\rho,-w_1w(\mu+\rho)\ra    \\
&=\ell\cdot\sum_{i}(w_1w(\mu+\rho))_i+\lan \la+\rho,-w_1 w(\mu+\rho)\ra. \\
\end{split}
\]

Now $t(\mu):=\sum_{i}(w_1w(\mu+\rho))_i=\sum_i(w(\mu+\rho))_i$
is an integer whose parity is the same as that of
$\sum_i(\mu+\rho)_i$ and depends only on $\mu$ (and the rank $k$),
and $q^\ell=\pm 1$ so $q^{\ell\cdot t(\mu)}=\pm1$ and we have
\[
\begin{split}
&\sum_{w\in W}\ve(w)q^{\lan\phi(\la)+\rho,w(\mu+\rho)\ra}
=\sum_{w\in W}\pm \ve(w) q^{\lan\la+\rho,-w_1w(\mu+\rho)\ra} \\
&=\pm \sum_{w^\prime \in W}
\ve(w^\prime)q^{\lan \la+\rho,w^\prime(\mu+\rho)\ra} \\
\end{split}
\]
where $w^\prime=-w_1w$.  Since the denominator of
$\dim_\F^{\mu}(V_\la)$ is independent of $\la$ the lemma is true
for $\mu \in P_+\cap\frac{1}{2}\Z^k\setminus \Z^k$. \hfill $\Box$

Let us pause for a moment to nail down exactly which sign
$\dim_\F^{\rho}(V_{\phi(\la)})=\dim_\F(V_{\phi(\la)})$ has in
terms of $\dim_\F^{\rho}(V_\la)$. Here there are two factors
governing signs of the characters: $\ve(-w_1)$ and the parity of
$\sum_iw(2\rho)_i$.  One has that:
$$\ve(-w_1)=\begin{cases}(-1)^{k/2}& \text{for $k$ even} \\
(-1)^{(k-1)/2} & \text{for $k$ odd}
\end{cases}$$

Furthermore we compute:
$q^{l\sum_iw(2\rho)_i}=(q^\ell)^k$
so we have the following result:
\begin{scholium}\label{dimschol} If $q^\ell=-1$ then
$$\dim_\F(V_{\phi(\la)})=
\begin{cases}\dim_\F(V_\la) &  k\equiv 0\bmod{4} \\
\dim_\F(V_\la) & k\equiv 1\bmod{4} \\
-\dim_\F(V_\la) & k\equiv 2\bmod{4} \\
-\dim_\F(V_\la) & k\equiv 3\bmod{4}
\end{cases}$$
Whereas if $q^\ell=1$:
$$\dim_\F(V_{\phi(\la)})=
\begin{cases}\dim_\F(V_\la) &  k\equiv 0\bmod{4} \\
-\dim_\F(V_\la) & k\equiv 1\bmod{4} \\
-\dim_\F(V_\la) & k\equiv 2\bmod{4} \\
\dim_\F(V_\la) & k\equiv 3\bmod{4}
\end{cases}$$
\end{scholium}

The following important lemma gives the decomposition rule for
tensoring with the object in $\F$
labeled by $\ga$.
\begin{lemma}\label{gam}
$V_{\ga}\otimes V_\mu = V_{\phi(\mu)}$ for all $\mu \in C_\ell$.
\end{lemma}
$Proof.$ By Lemmas \ref{phi} and \ref{pos} we know that
$\dim_\F^{\La_k}(V_{\ga})=\dim_\F^{\La_k}(\1)=1$ since
$\phi(\mathbf{0})=\ga$.  So
$$\dim_\F^{\La_k}(V_{\ga}\otimes V_\mu)=\dim_\F^{\La_k}(V_{\mu})=
\dim_\F^{\La_k}(V_{\phi(\mu)}).$$ Recall from \ref{antisym} that
$$\dim\Hom_\F(V_\ga\otimes V_\mu,V_\nu)=N_{\ga\mu}^{\nu}=
\sum_{W_\nu}\ve(w)m_{\ga\mu}^{w\cdot\nu}$$
where $W_\nu=\{w\in W_l: w\cdot\nu \in P_+\}$ and
$$m_{\ga\mu}^{w\cdot\nu}=\dim\Hom_{U_q\mathfrak{so}_{2k+1}}(V_{\ga}\otimes
V_\mu,V_{w\cdot\nu}).$$
Observe that the weight
$\phi(\mu)=\ga-w_1(\mu)$ is in $C_\ell$ and
$m_{\ga\mu}^{\phi(\mu)}=1$ (see
 \ref{tensorreps}).  The only way that $V_\phi(\mu)$
might fail to appear in the $\F$ decomposition is if $\phi(\mu)$
were equal to a reflection (under the dot action of $W_\ell$) of
$\ga+\kappa$ for some $\kappa \in P(\mu)$ (notice this also covers
weights in other Weyl chambers).  To see that this is impossible,
we use a geometric argument, although it is really nothing more
than an adaptation of the classical outer multiplicity formula.
First note that $\ga$ is a positive distance from all walls of
reflection under the dot action of $W_\ell$. Next observe that the
straight line segment from $\ga$ to $\ga+\kappa$ has Euclidean
length $|\kappa| \leq |\mu|$. So the reflected piecewise linear
path from $\ga$ to $w\cdot(\ga+\kappa)$ will not be straight, and
will have total length $|\kappa|$ as well.  Thus the straight line
segment from $\ga$ to $w\cdot(\ga+\kappa)$ must have length
strictly less than $|\mu|$, whereas the straight line segment from
$\ga$ to $\phi(\mu)$ has length $|\mu|$. So $V_\phi(\mu)$ appears
in the $\F$ decomposition of $V_{\ga}\otimes V_\mu$. But since
$\dim_\F^{\La_k}$ is positive on $C_\ell$ and
$$\dim_\F^{\La_k}(V_{\phi(\mu)})=\dim_\F^{\La_k}(V_{\ga}\otimes V_\mu)=
\sum_{\nu}N_{\ga\mu}^{\nu}\dim_\F^{\La_k}(V_{\nu})$$ it is clear
that $V_{\phi(\mu)}$ is the {\it only} object that appears in the
decomposition. \hfill $\Box$

\begin{rmk} This result can also be derived from \cite{LeTuraev},
Remark 3.9.  Le and Turaev studied symmetries in more general
settings for topological applications.  The involution $\phi$ also
appears in slightly more general setting in the paper \cite{Sawin}
by S. Sawin.
\end{rmk}
\subsection{The Family $\mF$ Summarized}
Let us collect together the important facts mentioned so far:
\begin{enumerate}
\item[1.] For a fixed $k$ and $\ell$ the corresponding family $\mF$ of
categories has a common Grothendieck semiring, denoted $Gr(\F)$.
\item[2.] $Gr(\F)$ has a unique positive character.
\item[3.] The involution $\phi$ preserves characters of $Gr(\F)$
up to a sign, and is induced by tensoring with $V_\ga$.
\item[4.] $Gr(\F)$ has at most $|C_\ell|$ distinct characters
each of which is a simultaneous eigenvector of the set $\M$ of
matrices.
\end{enumerate}

\section{Braid Group Representations}
In this section we analyze the representations $\C \B_n\rightarrow
\End_\F(\Vn)$ with an eye towards realizing these centralizer
algebras as (specializations of) quotients of $BMW$-algebras
$C_f(r,q)$ which we will define below.

Recall from Example \ref{rule2} that every object labeled by an
integer weight appears in an \emph{even} and an \emph{odd} tensor power of the
object $V_{\La_1}$.  Introduce the object $V=V_{\La_1}\otimes
V_\ga\cong V_{\phi(\La_1)}$.  By Lemma \ref{gam} we see that
$V$ is a generator for the category $\F$, since $V_\ga^{\otimes
2}\cong\1$.  We saw before that $V_{\La_k}$ was also a generator,
but $V$ has the advantage that $V^{\otimes 2}$ \emph{always}
decomposes as the direct sum of 3 simple objects regardless of the
rank $k$.  So we can take advantage of a computation of Tuba and
Wenzl (\cite{TuW2} proof of Lemma 6.3):
\begin{lemma}\label{dimseig}
Suppose $\OO$ is a semisimple ribbon category generated by a
simple object $X$ and $X^{\otimes 2}\cong \1\oplus Y\oplus Z$ with
$Y$ and $Z$ simple objects.  If the eigenvalues of $c_{X,X}$ are
$c_1, c_2$ and $c_3$ respectively on $Y, Z$ and $\1$ then:
$$\dim_\OO(X)=\pm\frac{c_3^2+c_1c_2-c_3(c_1+c_2)}{c_3(c_1^{-1}+c_2^{-1})}.$$
\end{lemma}
The proof of this lemma relies upon Lemma \ref{parttrace} and
the explicit computations in \cite{TuW1}.

\subsection{BMW-Algebras}
The algebras $C(r,q)$ are
quotients of the group algebra of Artin's braid group $\B_f$ and
were studied extensively in \cite{BCD} and \cite{TuraevWenzl2}, and
more recently in \cite{TuW2}.
\begin{definition}\label{defrels}
Let $r,q \in \C$ and $f \in \N$, then $C_f(r,q)$ is the
$\C$-algebra with invertible generators $g_1,g_2,\ldots,g_{f-1}$
and relations:
\begin{enumerate}
\item[(B1)] $g_ig_{i+1}g_i=g_{i+1}g_ig_{i+1}$,
\item[(B2)] $g_ig_j=g_jg_i$ if $|i-j|\geq 2$,
\item[(R1)]
$e_ig_i=r^{-1}e_i$,
\item[(R2)] $e_ig_{i-1}^{\pm 1}e_i=r^{\pm 1}e_i$,

where $e_i$ is defined by
\item[(E1)] $(q-q^{-1})(1-e_i)=g_i-g_i^{-1}$
\end{enumerate}
\end{definition}
Notice that (E1) and (R1) imply
$$(g_i-r^{-1})(g_i-q)(g_i+q^{-1})=0$$ for all $i$.  So the image
of $g_i$ on any finite dimensional representation has at most three
eigenvalues, which are distinct if $q^2\not=-1$ and $r\not=\pm q^{\pm 1}$.
Notice further that the image of
$e_i$ is a multiple of the projection onto the
$g_i$-eigenspace corresponding to eigenvalue $r^{-1}$.

There
exists a trace $tr$ on the family of algebras $C_f(r,q)$ uniquely
determined by the values on the generators, and inductively
defined by (see \cite{BCD}):
\begin{enumerate}
\item[(T1)]  $tr(1)=1$
\item[(T2)] $tr(g_i)=r(\frac{q-q^{-1}}{r-r^{-1}+q-q^{-1}})$,
\item[(T3)] $tr(axb)=tr(ab)tr(x)$ for $a,b \in C_{f-1}(r,q)$ and $x \in
\{g_{f-1},e_{f-1},1\}$.
\end{enumerate}

The existence of such a trace comes from
the well-known Kauffmann link invariant.  When $q$ is a root of
unity and $r$ is plus or minus a power of $q$ then the algebras
$E_f:=C_f(r,q)/Ann(tr)$ are finite dimensional and semisimple and
hence isomorphic to a direct sum of full matrix algebras.

In \cite{TuraevWenzl2} the authors construct a family of self-dual
Hermitian ribbon categories from the sequence of algebras
$\C\subset\cdots E_f\subset E_{f+1}\subset\cdots$ for various
choices of $r$ and $q$.  The objects in these categories are the
idempotents in the algebras $E_f$, $f\geq 1$ and the morphisms are
images of \emph{tangles}.  Since the algebra $E_f$ is a quotient
of the group algebra of the braid group $\C B_f$, the braiding in
the category is obtained directly as images of elements in $E_f$.
The construction is quite involved, so we will be content to
outline the important properties leaving the interested reader to
seek details in the above reference as well as \cite{TuW2}.

\subsection{The $BC$-Case and the Family $\mV$}\label{Vcat} Fix
$q$ with $q^2$ a primitive $\ell$th root of unity, and let
$r=-q^{2k}$.  Denote by $\V$ the corresponding self-dual Hermitian
ribbon category as constructed in \cite{TuraevWenzl2}. This is
known as the \emph{ortho-symplectic} or $BC$-case in the
literature.  For $\ell$ odd, we have the following:
\begin{enumerate}
\item[1.]  The simple objects of $\V$ are labeled by Ferrer's
diagrams
$\la\in\Ga(k,\ell)$ where:
$$\Ga(k,\ell):=\{\la:\la_1^{\prime}+\la_2^\prime \leq 2k+1,\la_1\leq(\ell-2k-1)/2\}$$
with $\la_i$ (resp. $\la_i^\prime$) the number of boxes in the
$i$th row (resp. column) of $\la$ (see \cite{BCD}). \item[2.] The
object $X:=X_\Box$ generates $\V$ (see \cite{BCD}). \item[3.] For
$\mu\in\Ga(k,\ell)$, $X_\mu$ is a simple subobject of the tensor
product $X\otimes X_\la$ if and only if $\mu$ can be obtained by
adding/deleting one box to/from $\la$ (see \cite{TuraevWenzl2}).
\item[4.] $Gr(\V)\cong Gr(Rep(O(2k+1)))/\mathcal{J}$ where
$\mathcal{J}$ is some ideal.  (see \cite{TuW2}). \item[5.]
$\dim_\V(X)=\frac{[-2k]}{[1]}+1$  (see \cite{BCD}). \item[6.] The
eigenvalues of the braiding morphism $c_{X,X}$ on the simple
subobjects $\{\one, X_{[2]}, X_{[1^2]}\}$ are respectively either
$\{-q^{-2k},q,-q^{-1}\}$ or $\{-iq^{-2k}, iq, -iq^{-1}\}$.
(depending on a choice of a braiding, see \cite{TuW2}). \item[7.]
By replacing the braiding morphism $c_{X,X}$ by its negative,
inverse or negative-inverse we get 3 new inequivalent ribbon
categories with the same Grothendieck semiring as $\V$. \item[8.]
There is an algebra isomorphism $\End_\V(X^{\otimes n})\cong E_n$
that preserves the $\C \B_n$-module structure (see \cite{TuW2}).
\end{enumerate}
The key theorem we will use is the following
 special case of the main result in \cite{TuW2} (Theorem 9.5):
\begin{prop}\label{mainTuW2}
Fix $k$ and $\ell$.  Let $\mV$ be the family of braided tensor
categories constructed from $BMW$-algebras where $q^2$ is any
primitive $\ell$th root of unity and $r=-q^{2k}$, and $c_{X,X}$ is
one of the four braiding morphisms as in item 7 above.  Then any
braided tensor category $\OO$ with $Gr(\OO)\cong Gr(\V)$ such that
the braiding morphism $c_{Y,Y}$ where $Y\in\OO$ is the object
corresponding to $X\in\V$ has 3 distinct eigenvalues is equivalent
(as a braided tensor category) to a member of the family $\mV$.
\end{prop}

\section{Main Theorem}
We now proceed to prove:
\begin{theorem}\label{maintheorem}
Fix $\F\in \mF$. Then as a braided tensor category, $\F$ is
equivalent to some $\V\in\mV$.
\end{theorem}
The proof is outlined as follows:
\begin{enumerate}
\item[Step 0] It is sufficient to show the theorem for any
fixed $\F\in\mF$.
\item[Step 1] The image of $\C \B_n$ in $\End_\F(\Vn)$ is
a quotient of $E_n$.
\item[Step 2] There exists a $\V\in\mV$ such that as algebras
$\End_\V(X^{\otimes n})\cong\End_\F(\Vn)$.
\item[Step 3] $Gr(\V)\cong Gr(\F)$.
\item[Step 4]  As braided tensor categories $\V$ and $\F$ are
equivalent up to 4 possible choices of braiding morphism
$c_{V,V}$.
\end{enumerate}
\begin{rmk}
It should be emphasized that the fusion rules of $\mF$ are \emph{a
priori} only obtained as a quotient of the representation category
of $O(2k+1)$ which is the ``integer half" of the representation
category of $\mathfrak{so}_{2k+1}$ (i.e. lacking the spinnor
representations).  So while it was well-known that there is a
relationship between $\mV$ and $\mF$ (see \cite{BCD}), we show
that \emph{all} of $\mF$ can be obtained as a quotient of the
Turaev-Wenzl category.
\end{rmk}

\subsection{Proof of Step 0} Since all $\F\in\mF$ share the
same Grothendieck semiring and the eigenvalues are distinct by
Proposition \ref{reigs} (we explicitly compute them below),
\ref{mainTuW2} implies that once we have established the theorem
for some $\F\in\mF$ we will be done.
\subsection{Proof of Step
1} Step 1 will follow as soon as we show that the images of the
braid generators: $R_i:=\one_V^{\otimes (i-1)}\otimes
c_{V,V}\otimes\one_V^{\otimes (n-i-1)}$ satisfy the defining
relations of $C_n(r,q)$ as well as the trace conditions (see
\ref{defrels}) for appropriate choice of $r$ and $q$.

The object $V^{\otimes 2}$ decomposes as the sum of the three
objects $\1$, $V_1:=V_{(2,0,\ldots,0)}$ and
$V_2:=V_{(1,1,0,\ldots,0)}$. Applying Proposition \ref{reigs} we
see that the eigenvalues of $(c_{V,V})^2$ on $V_1$, $V_2$ and $\1$
depend on the parity of $k$ and the sign of $q^\ell=\pm 1$ as
follows:
\begin{enumerate}
\item[1.] On $V_1$: \begin{equation}
c_1^2=
\begin{cases}-q^{-4}& \text{if $k$ odd and $q^\ell=-1$} \\
q^{-4} & \text{otherwise}
\end{cases}
\end{equation}
\item[2.] On $V_2$: \begin{equation}
c_2^2=
\begin{cases}-q^{4}& \text{if $k$ odd and $q^\ell=-1$} \\
q^{4} & \text{otherwise}
\end{cases}
\end{equation}

\item[3.] On $\1$: \begin{equation}
c_3^2=
\begin{cases}-q^{-8k}& \text{if $k$ odd and $q^\ell=-1$} \\
q^{-8k} & \text{otherwise}
\end{cases}
\end{equation}
\end{enumerate}
So the eigenvalues are $\{c_1,c_2,c_3\}$ either $\{\pm q^{-2},\pm
q^2,\pm q^{-4k}\}$ or $\{\pm iq^{-2},\pm iq^2,\pm iq^{-4k}\}$
where the sign choices are independent.  For simplicity (by Step
0) we assume that $q^\ell=-1$.

 By Lemma \ref{phi} we have that
$\dim_\F(V)=\pm\dim_\F(V_{\La_1})$ and an easy computation using
the equation for $\dim_\F$ in \ref{qdim} we get:
\begin{equation}\label{Vdim}
\pm\dim_\F(V)=\frac{[4k]}{[2]}+1
\end{equation}

Next we make a change of parameter: $\tq\rightarrow -q^2$. Observe
that $\tq$ is still a primitive $\ell$th root of unity with
$\tq^\ell=-1$, so by Step 0 we can proceed with this altered
category.  This change gives us:
\begin{equation}
\pm\dim_\F(V)=\frac{-[2k]_{\tq}}{[1]_{\tq}}+1
\end{equation}
and $\{c_1,c_2,c_3\}$ is either $\{\pm \tq^{-1},\pm \tq,\pm
\tq^{-2k}\}$ or $\{\pm i\tq^{-1},\pm i\tq,\pm i\tq^{-2k}\}$.

Using Lemma \ref{dimseig} we can test the possible choices of
$\{c_1,c_2,c_3\}$ by computing $\pm\dim_\F(V)$ and comparing with
$\frac{-[2k]_{\tq}}{[1]_{\tq}}+1$.  A (somewhat tedious)
computation forces $\{c_1,c_2,c_3\}$ to be one of the two choices
$\pm\{ -\tq^{-1},\tq,-\tq^{-2k}\}$ for $k$ even and
$\pm\{-i\tq^{-1},i\tq,-i\tq^{-2k}\}$ for $k$ odd.  Now by changing
the sign of $c_{X,X}$, we can change the sign of the corresponding
eigenvalues for the target category $\V$, so we assume the
eigenvalues $\{c_1,c_2,c_3\}$ are $\{ -\tq^{-1},\tq,-\tq^{-2k}\}$
for $k$ even and $\{-i\tq^{-1},i\tq,-i\tq^{-2k}\}$ for $k$ odd. So
comparing with the defining relations for the $BMW$-algebras and
setting $-\tq^{2k}=r$ we need to show (for the $k$ even case):
\begin{enumerate}
\item[(B1)] $R_iR_{i+1}R_i=R_{i+1}R_iR_{i+1}$,
\item[(B2)] $R_iR_j=R_jR_i$ if $|i-j|\geq 2$,
\item[(R1)]$S_iR_i=r^{-1}S_i$,
\item[(R2)] $S_iR_{i-1}^{\pm 1}S_i=r^{\pm 1}S_i$,

where $S_i$ is defined by
\item[(E1)] $(\tq-\tq^{-1})(\one-S_i)=R_i-R_i^{-1}$
\item[(T1)]  $tr_\F(\one_V)=1$
\item[(T2)] $tr_\F(R_i)=r(\frac{\tq-\tq^{-1}}{r-r^{-1}+\tq-\tq^{-1}})$,
\item[(T3)] $
tr_\F((a\otimes \one_V)x(b\otimes \one_V))= tr_\F((a\circ
b)\otimes \one_V)tr_\F(x)$ for $a,b \in \End_\F(V^{\otimes
(f-1)})$ and $x \in \{R_{f-1},S_{f-1},1\}\subset\End(V^{\otimes
f})$.
\end{enumerate}

Here $tr_\F(a):=Tr_\F(a)/\dim_\F(\Vn)$ where $a\in\End_\F(\Vn)$.
Relations (B1) and (B2) are immediate from the braiding axioms,
and (T1) follows from the definition of the normalized trace
$tr_\F$.  Relation (R1) follows from the computation of the
eigenvalues of $R_i$ and definition (E1).  To verify (R2) it is
sufficient to consider $i=2$ and verify the relation on
$\End_\F(V^{\otimes 3})$.  Since $S_1\in\End_\F(V^{\otimes 2})$ is
$1+\frac{r-r^{-1}}{\tq-\tq^{-1}}$ times the projection onto the
subobject $\1$ in $V^{\otimes 2}$, we can apply Lemma
\ref{parttrace} to get:
$$S_2R_1S_2=tr_\F(c_{V,V})(1+\frac{r-r^{-1}}{\tq-\tq^{-1}})S_2$$
so (R2) will follow from (T2).  Applying \ref{parttracecor} to the
eigenspace decomposition $R_1=r^{-1}p_{\1}-\tq^{-1}p_{V_1}+\tq
p_{V_2}$ acting on $V\otimes V$ and computing $\dim_\F(V_1)$ and
$\dim_\F(V_2)$ from the definition it is a matter of simple
algebra to verify (T2).  All that remains is to verify (T3).  But
since the algebras $\End(V^{\otimes f})$ are semisimple and finite
dimensional it is enough to show (T3) for $a,b$ minimal
idempotents.  But this reduces to Lemma \ref{markov}.

So we conclude that $\End(\Vn)$ contains a quotient of $E_n$ for
all $n$.\hfill $\Box$
\subsection{Proof of Step 2} Using the
fact that $E_n\cong\End_\V(X^{\otimes n})$ for all $n$ together with Step
1 we need only show that $\dim\End_\F(\Vn)=\dim\End_\V(X^{\otimes
n})$ for \emph{any} $\V$ and $\F$ in their respective families to
conclude that the action of $\C \B_n$ on these two algebras is the
same.

Several tensor categories
will be bandied about in what follows.
Recall first the following sets:
\begin{enumerate}
\item $\Gamma(k,\ell)=
\{\la:\la_1^{\prime}+\la_2^{\prime}\leq
2k+1,\la_1\leq(\ell-2k-1)/2\}$.  Here $\la$ is a Ferrer's diagram,
and $\la_i^\prime$ is the number of boxes in the $i$th column.
\item $P_+=
\{\la\in\Z^k\cup(\Z^k+\frac{1}{2}(1,1,\ldots,1)):\la_1\geq\la_2\geq\ldots\la_k\geq 0\}$
\item $C_\ell=\{\la\in P_+: \frac{\ell-2k}{2}\geq\la_1\}$.
\end{enumerate}
Table 1 will serve as a lexicon of notation.  The first column is
the category, the second the labeling set for simple objects, and
the third the notation used for the simple object labeled by
$\la$.
\begin{table}\label{cattable}
\caption{Tensor Categories}
\begin{tabular}{*{3}{|c}|}
\hline
Category & Labeling Set & Objects \\
\hline\hline
$Rep(O(2k+1))$ &  Diagrams $\la$, $\la^{\prime}_1+\la^\prime_2\leq 2k+1$ & $W_\la$ \\
\hline
$\V$ & $\Ga(k,\ell)$  & $X_\la$ \\
\hline
$Rep(U_q\mathfrak{so}_{2k+1})$, $|q|\not=1$ & $P_+$ & $V_\la$ \\
\hline
$\F$ & $C_\ell$ & $V_\la$\\
\hline
\end{tabular}
\end{table}

Next we note a few homomorphisms that exist between the Grothendieck
semirings of these tensor categories.
\begin{enumerate}
\item[1.] As we mentioned \ref{Vcat}, the ring $Gr(\V)$ is a
quotient of $Gr(Rep(O(2k+1)))$. Provided $\mu\in\Ga(k,\ell)$ we
have:
$$\dim\Hom_\V(X_\mu,X_{\Box}\otimes X_\la)=
\dim\Hom_{O(2k+1)}(W_\mu,W_{\Box}\otimes W_\la).$$

\item[2.] Define a map from the set of $O(2k+1)$ dominant weights
(Ferrer's diagrams with at most $2k+1$ boxes in the first two
columns) to the integer weights of $\mathfrak{so}_{2k+1}$ by
restricting and differentiating the irreducible representations.
Explicitly this associates to $\la$ the Ferrer's diagram
$\overline{\la}$ identical to $\la$ except the first column has
$\min\{2k+1-\la_1^\prime,\la_1^\prime\}$ boxes (here
$\la_1^\prime$ is the number of boxes in the first column of
$\la$).  By filling in zeros for empty rows, we express
$\overline{\la}$ as a $k$-tuple in our standard notation for
dominant weights of $\mathfrak{so}_{2k+1}$.  The map
$\la\rightarrow \overline{\la}$ induces a homomorphism from
$Gr(Rep(O(2k+1)))$ to $Gr(Rep(\mathfrak{so}_{2k+1}))$. From this
we deduce:
$$\dim\Hom_\V(W_\mu,W_{\Box}\otimes W_\la)=
\dim\Hom_{\mathfrak{so}_{2k+1}}(V_{\overline{\mu}},V_{\La_1}\otimes
V_{\overline{\la}}).$$

\item[3.] For generic $q$, the semirings
$Gr(Rep(\mathfrak{so}_{2k+1}))$ and
$Gr(Rep(U_q\mathfrak{so}_{2k+1}))$ are isomorphic. For this reason
we denote the simple objects from both categories by $V_\la$.

\item[4.] The category $\F$ is obtained from
$Rep(U_q\mathfrak{so}_{2k+1})$ as a quotient.  Heedless of any
potential confusion, we denote the simple objects in $\F$ by
$V_\la$ as well. Recall from Example \ref{rule2} that for any
integer weight $\la\in C_\ell$:
$$V_\mu\subset V_{\La_1}\otimes V_\la \iff
V_\mu\subset V_{\La_1}\otimes V_\la,\ \mu \in C_\ell.$$
\end{enumerate}

Define a bijection $\Psi: \Gamma(k,\ell) \rightarrow C_\ell$ by
\begin{equation}
\Psi(\la)=
\begin{cases}
\overline{\la},   &\text{if $|\la|$ is even,} \\
\phi(\overline{\la}),  &\text{if $|\la|$ is odd.}
\end{cases}
\end{equation}
Observing that the tensor product of any simple object in $\V$ (resp. $\F$)
with the generating object $X$ (resp. $V$) is multiplicity free,
the algebras $\End_\V(X^{\otimes n})$ and $\End_\F(\Vn)$
are isomorphic once we prove:
\begin{lemma}\label{homeq}
Let $\mu,\la\in\Ga(k,\ell)$.  Then
$$\dim\Hom_\V(X_\mu,X\otimes X_\la)=
\dim\Hom_\F(V_{\Psi(\mu)},V\otimes V_{\Psi(\la)}).$$
\end{lemma}
$Proof.$ Using the first homomorphism of Grothendieck semirings
above and the assumption that $\mu\in\Ga(k,\ell)$, we have
$$\dim\Hom_\V(X_\mu,X\otimes X_\la)=
\dim\Hom_{O(2k+1)}(W_\mu,W_{\Box}\otimes W_\la).$$ Restricting to
$SO(2k+1)$, differentiating and applying the third homomorphism
above we have
$$\dim\Hom_{U_q\mathfrak{so}_{2k+1}}
(V_{\overline{\mu}},V_{\La_1}\otimes V_{\overline{\la}})=
\dim\Hom_{O(2k+1)}(W_\mu,W_{\Box}\otimes W_\la).$$
Now we split into the two cases from the definition of $\Psi$:\\
Case I: $|\la|$ is even (so $|\mu|$ is odd)\\
Since $\overline{\mu}\in C_\ell$ and $V_{\Psi(\la)}=V_{\overline{\la}}$
we see that \begin{equation}\label{laeven}
\dim\Hom_\F(V_{\overline{\mu}},V_{\La_1}\otimes V_{\Psi(\la)})=
\dim\Hom_{U_q\mathfrak{so}_{2k+1}}
(V_{\overline{\mu}},V_{\La_1}\otimes V_{\overline{\la}})
\end{equation}
Lemma \ref{gam} implies that
$V_\ga\otimes V_{\overline{\mu}}=V_{\phi(\overline{\mu})}=V_{\Psi(\mu)}$
as objects in $\F$, and similarly $V_\ga\otimes V_{La_1}=V$.
So tensoring with
$V_\ga$ (see example \ref{laeven}) we have:
$$\dim\Hom_\F(V_{\Psi(\mu)},V\otimes V_{\Psi(\la)})=
\dim\Hom_\F(V_{\overline{\mu}},V_{\La_1}\otimes V_{\Psi(\la)}).$$
Case II: $|\la|$ is odd (so $|\mu|$ is even)\\
In this case $V_{\Psi(\la)}=V_\ga\otimes V_{\overline{\la}}$
and $V_{\Psi(\mu)}=V_{\overline{\mu}}$ so using the fact that
$V_\ga\otimes V_\ga=\1$ we derive similarly that
$$\dim\Hom_\F(V_{\Psi(\mu)},V\otimes V_{\Psi(\la)})=
\dim\Hom_\F(V_{\overline{\mu}},V_{\La_1}\otimes V_{\Psi(\la)})$$
in this case.  \hfill $\Box$\\

This lemma implies that $$\dim\Hom_\V(X_\nu,X^{\otimes n})=
\dim\Hom_\F(V_{\Psi(\nu)},\Vn)$$ by an easy induction argument.
Thus we have shown that $\dim\End_\F(\Vn)=\dim\End_\V(X^{\otimes
n})$.\hfill $\Box$
\subsection{Proof of Step 3} In Step 1 we
established that the action of $\C \B_n$ on $\End_\F(\Vn)$ factors
through $E_n$ which is isomorphic to $\End_\V(X^{\otimes n})$ both
as algebras and $\C \B_n$-modules. Combining this with Step 2 we
conclude that $\End_\F(\Vn)$ and $\End_\V(X^{\otimes n})$ are
isomorphic both as algebras and as $\C \B_n$-modules. Since $V$
and $X$ generate their respective categories, this implies that
the Grothendieck semirings $Gr(\F)$ and $Gr(\V)$ are isomorphic.
In fact, it is tedious but straightforward to show that:
\begin{cor}
$\Psi$ induces an isomorphism $Gr(\V)\cong Gr(\F)$.
\end{cor}
That is, the map $\Psi$ defined above on the labeling sets of
simple objects describes precisely the correspondence between
these two categories. \hfill $\Box$
Observe that we also get the following theorem as a consequence
(see \cite{OW} for similar statements):
\begin{theorem}
The centralizer algebra $\End_\F(\Vn)$ is generated by the image
of $\C \B_n$.
\end{theorem}
\subsection{Step 4} Since any $\F\in\mF$ has the same
Groethendieck semiring as any $\V\in\mV$ and the braiding morphism
$c_{V,V}$ has 3 distinct eigenvalues, considered as braided tensor
categories, the family $\mF$ is a subfamily of $\mV$ by
Proposition \ref{mainTuW2}. \hfill $\Box$

\subsection{Extension to Lie type $C$}
 It is known that the Turaev-Wenzl categories of
type $BC$ have the same Grothendieck semiring as the categories
corresponding to quantum groups of Lie type $C$ at odd roots of
unity (see \cite{BB}).  Combining this with our result, we get the
following rank-level duality type corollary:
\begin{cor}\label{ranklevel}
 The ribbon
categories corresponding to the rank $k$ quantum group of Lie type
$B$ and the rank $(\ell-2k-1)/2$ quantum group of Lie type $C$ at
a $\ell$th root of unity have the same tensor product rules.
\end{cor}
Moreover, we can compute the eigenvalues of the braiding
isomorphism $c_{V,V}$ for $V$ the highest weight quantum group
module of type $C_r$ corresponding to the weight $(1,0,\ldots,0)$.
Here it is even easier than for type $B$ as we can use
\cite{ramleduc} Corollary 2.22(3).  The eigenvalues are:
$$\{q,-q^{-1},-q^{-2r-1}\}.$$
Using Corollary \ref{ranklevel} we set $r=(\ell-2k-1)/2$ which
gives us eigenvalues $$\{q,-q^{-1},-q^\ell q^{-2k}\}$$ which can
be made to match those of $\mV$ by changing an overall sign and/or
transposing all Young diagrams as in \cite{TuW2}.  Thus we can
apply the theorem of Tuba and Wenzl to see that the Lie type $C$
at odd roots of unity categories can be included in this family of
ribbon categories.
\section{Failure of Unitarity}
We will show that no member of the family of braided tensor
categories $\mV$ can have the structure of a Unitary ribbon
category. We showed in Lemma \ref{Unique} that there is a unique
positive character for the Grothendieck semiring $Gr(\F)$.  By the
above equivalence, we also know that the same true for $Gr(\V)$
for any $\V\in\mV$.  Lemma \ref{dimseig} shows that $\dim_\V(X)$
is uniquely determined up to a sign by the eigenvalues of the
braiding morphism $c_{X,X}$ and so we have that:
\begin{equation}
\pm\dim_\V(X)=\frac{-[2k]}{[1]}+1
\end{equation}
So if we can show that $\pm\dim_\V(X)$ is
never equal to the unique positive character of Lemma
\ref{Unique} above for any choice of $q^2$ a primitive $\ell$th root of unity then
we will have shown that this abstract category does not support
both positivity and a braiding.
For any $\la\in\Ga(k,\ell)$ we denote by $\Dim(\la)$ the unique positive character of
$Gr(\V)$. Furthermore, we set
$$f_{\la}(z)=\dim_\V(X_\la)\mid_{(e^{z\pi i/\ell})}$$
for $1\leq z \leq \ell-1$ and $\gcd(\ell,z)=1$ so that
$f_{\la}(z)$ takes on all possible values of $\dim_\V(X_\la)$ as
$\V$ ranges over the family $\mV$.
 We may now formulate:
\begin{theorem}
If $2(2k+1)<\ell$ then $f_{\la}(z)\not=\Dim(\la)$ for any $z$ with
$1\leq z \leq \ell-1$ and $\gcd(\ell,z)=1$.
\end{theorem}
  Since both $f_{\la}(z)$ and $\Dim(\la)$ are both characters of
  $Gr(\V)$ (i.e. they are
normalized so that their values at the trivial object are 1),
this theorem will be a consequence of the following:
\begin{lemma}
Let $h(z)=f_{\Box}(z)$.  Then if $2(2k+1) \leq \ell$ and $1 \leq z
\leq \ell-1$ with $\gcd(\ell,z)=1$ then $|h(z)|<\Dim(\Box)$.
\end{lemma}
$Proof.$  We start by showing that $h(z)<\textrm{Dim}_q(\Box)$.
We have that $h(z)=\frac{-\sin(2 k z \pi/\ell)}{\sin(z\pi/\ell))} +1$ and
$\textrm{Dim}_q(\Box)=\frac{\sin((2k+1)\pi/\ell)}{\sin(\pi/\ell)}$.  First one notes
that $\textrm{Dim}_q(\Box) >1$ and so $h(z) < \textrm{Dim}_q(\Box)$ if $z \leq \ell/2k$.
So the lemma is true for $z \in I_1=[1,\ell/2k]$.

Next we make a change of variables $z \rightarrow \ell-z^{\prime}$
in order to eliminate large $z$.  We define
$$g(z^{\prime})=h(\ell-z^{\prime})=
\frac{\sin(2kz^{\prime}\pi/\ell)}{\sin(z^{\prime}\pi/\ell)}+1$$
with $1 \leq z^\prime \leq \ell-1$.  Using the sum expansion of
$\frac{q^{2k}-q^{-2k}}{q-q^{-1}}$ we can write
$$g(z^{\prime})=1+2\sum_{1\leq j\leq k}\cos((2j-1)z^{\prime}\pi/\ell).$$
By taking a derivative of $g(z^{\prime})$ we find that it is a decreasing
function of $z^\prime$ on the interval $I_2^\prime=[2,\frac{\ell}{2k-1}]$,
which is nonempty if $2(2k-1) \leq \ell$.  Thus if
$g(2) < \Dim(\Box)$ then $g(z^\prime) < \Dim(\Box)$ on all of $I_2^\prime$.
Expanding $\Dim(\Box)$ we compute:
\[
\Dim(\Box)-g(2)=
2\sum_{1\leq j \leq k}[\cos(2j\pi/\ell)-\cos(2(2j-1)\pi/\ell)].
\]
Using the trigonometric formulas found in the back of any calculus book
we may express each term $\cos(2j\pi/\ell)-\cos(2(2j-1)\pi/\ell)$
as $2\sin((3j-1)\pi/\ell)\sin((j-1)\pi/\ell)$.  Provided
$3j-1 \leq 3k-1 \leq \ell$, each of these terms is positive.
But we already have the stronger restriction $2(2k+1) \leq \ell$, thus
we have $g(z^\prime) < \Dim(\Box)$ on $I_2^\prime$ that is,
$h(z) < \Dim(\Box)$ on $I_2=[\ell-\frac{\ell}{2k-1},\ell-2]$.
We check the case $z^\prime=1$ separately:
$$\Dim(\Box)-g(1)=\sum_{1\leq j\leq k}[\cos(2j\pi/\ell)-\cos((2j-1)\pi/\ell]$$
and each term can be factored as:
$$-2\sin(\pi/2\ell)\sin((4j-1)\pi/2\ell)$$
which is always strictly negative since $4j-1<2\ell$ for all $j\leq k$.

The only remaining $z$ to eliminate are those between $I_1$ and $I_2$.
To this end we use the following estimates which come from approximating
$\sin(x)$ from below by $1-|2x/\pi-1|$ on the interval $0\leq x \leq \pi$:
$$h(z)<\frac{1}{\sin(z\pi/\ell)}+1<2(2k+1)/\pi\leq \Dim(\Box)$$
which are valid for
$z \in I_3=[\frac{\ell\pi}{4(2k+1)-2\pi},\ell-\frac{\ell\pi}{4(2k+1)-2\pi}]$
provided $2(2k+1) < \ell$.  It is now easy to see that $[1,\ell-2] \cup \{\ell-1\}\subset
I_1\cup I_2\cup I_3$ thus proving that $h(z)<\Dim(\Box)$ for any
$z,\ell,k$ as in the statement.

With a few modifications to this proof we can show that
$-h(z)<\textrm{Dim}_q(\Box)$ as follows.  On $I_3$ our estimates are still
valid.  We observe that $-h(z)$ is decreasing on
$[1,\frac{\ell}{2k-1}]$ so one need only check that
$\Dim(\Box)>-h(1)$, which is straightforward.
By changing variables as we did above we can also eliminate $z \in
[\ell-\frac{\ell}{2k},\ell-2]$ using the observation that
$\Dim(\Box) >1$ again.  One must again check the case $z=\ell-1$ separately but the
same basic argument works as above except we must use the stronger condition
$4k-1\leq \ell$ since the factors involved are cosines.\hfill $\Box$\\

So we have shown that for no $q^2$ a primitive $\ell$th root of unity does the
categorical dimension of any $\V\in\mV$ achieve the value of the unique positive
character of $Gr(\V)$ (or $Gr(\F)$).  Observe that
In order to apply Theorem \ref{nounitary} and complete the proof, we
observe that there is a simple object $X_\tau\in\V$ with $|\tau|$ even and
$\dim_\V(X_\tau)<0$.
For if all simple $X_\tau$ with $|\tau|$ even had positive
dimension by multiplying by $(-1)^{|\la|}$ we would get a
$\dim_\V$ function that was positive on all simple objects but
with the same $\dim_\V(X_\Box)$ up to a sign, which is impossible
by the above lemma.  Since every $X_\tau$ with $|\tau|$ even
appears in an even power of $X$, we can apply Theorem
\ref{nounitary} and conclude that:
\begin{cor}
No braided tensor category $\OO$ with $Gr(\OO)\cong Gr(\V)$ (or
$Gr(\F)$) is unitarizable.
\end{cor}
\begin{rmk}
  It should be noted that it
   was previously thought that the Turaev-Wenzl categories in
the $BC$-case \emph{are} unitary for the choices $q=-e^{\pm\pi
i/\ell}$ (see \cite{TuraevWenzl2}, Theorem 11.2). The critical
theorem used to prove the positivity of the form is in \cite{BCD},
Theorem 6.4. However, the discovery of a slight miscalculation in
the case $\ell$ odd reveals that the argument fails in the present
case.
\end{rmk}

\bibliographystyle{abbrv}

\begin{thebibliography}{AA}
\bibitem[A]{andersen} H.H. Andersen, \emph{Tensor products of quantized
tilting modules}, Comm. Math. Phys. \textbf{149} (1991), 149-159.
\bibitem[AP]{AndPar} H.H. Andersen and J. Paradowski,
\emph{Fusion categories arising from semisimple Lie algebras},
Comm. Math. Phys. \textbf{169} (1995) 563-588.
\bibitem[BK]{BK} B. Bakalov and A. Kirillov Jr., Lectures on tensor categories
and modular functors. AMS Providence, 2001.
\bibitem[BB]{BB} A. Beliakova and C. Blanchet, \emph{Modular categories of
type B, C and D}, Comment. Math. Helv. \textbf{76} (2001) 467-500.
\bibitem[BW]{BW} J. Birman and H. Wenzl, \emph{Braids, link polynomials and a new algebra},
Trans. AMS \textbf{313} (1989) 249-273.
\bibitem[ChPr]{ChPr} V. Chari and A. Pressley, A Guide to Quantum Groups. Cambridge
University Press, Cambridge, 1994.
\bibitem[D]{D} V. G. Drinfeld, \emph{Almost cocommutative Hopf
algebras}, Leningrad Math. J.  \textbf{1} (1990), no. 2, 321--342
(translation).

\bibitem[Ga]{Ga} F. R. Gantmacher, The Theory of Matrices, v. 2, Chelsea, New
York, 1959.

\bibitem[GWa]{bigred} R. Goodman and N. Wallach, Representations and Invariants
of the Classical Groups, Cambridge University Press, Cambridge, 1998.
\bibitem[Hu]{humphreys} J. E. Humphreys, Introduction to Lie Algebras and
Representation Theory, Springer-Verlag, New York.
\bibitem[Ja]{jantzen} J. C. Jantzen, Lectures on Quantum Groups, AMS 1996.
\bibitem[K]{K}  C. Kassel, Quantum Groups. Graduate Texts in Mathematics, 155.
Springer-Verlag, New York, 1995.
\bibitem[Kf]{Kauf} L. H. Kauffmann, \emph{An invariant of regular isotopy},
Trans. Amer. Math. Soc. \textbf{318} (1990) no. 2, 417-471.
\bibitem[Ki]{kirillov} A. Kirillov Jr., \emph{On an inner product in modular categories},
J. of AMS \textbf{9} (1996), 1135-1170.
\bibitem[LT]{LeTuraev} T. Le and V. Turaev, \emph{Quantum groups and
ribbon $G$-categories}, J. Pure Appl. Algebra, \textbf{178} (2003)
169-185.
\bibitem[Lu]{lusztig} G. Lusztig, Introduction to Quantum Groups, Birkh\"auser,
1993.
\bibitem[M]{M} J. Murakami, \emph{The Kaufmann polynomial of links and representation theory},
Osaka J. Math. \textbf{24} (1987), 745-758.
\bibitem[OW]{OW} R. C. Orellana and H. Wenzl, \emph{$q$-centralizer algebras
for spin groups}, J. Algebra, \textbf{253} (2) (2002) 237-275.
\bibitem[LR]{ramleduc} R. Leduc and A. Ram, A\emph{ ribbon Hopf algebra
approach to the irreducible representations of centralizer
algebras: The Brauer, Birman-Wenzl and type $A$ Iwahori-Hecke
algebras}, Adv. Math. \textbf{125} (1997) 1-94.
\bibitem[R]{thesis} E. Rowell, Tensor categories arising from
quantum groups and BMW-algebras at odd roots of unity, thesis,
U.C. San Diego, 2003.
 \bibitem[S]{Sawin} S. Sawin,
 \emph{Jones-Witten invariants for nonsimply connected Lie groups and the geometry
  of the Weyl alcove},
Adv. Math. \textbf{165} (2002), no. 1, 1-34.
\bibitem[S2]{S2} S. Sawin, Quantum groups at roots of unity and modularity, preprint, arXiv: math.QA/0308281.
\bibitem[Tur]{Tur} V. G. Turaev, Quantum Invariants of Knots and 3-Manifolds.
 de Gruyter Studies in Mathematics, 18. Walter de Gruyter \& Co., Berlin, 1994.
\bibitem[TW1]{TuraevWenzl1} V. G. Turaev, H. Wenzl, \emph{Quantum invariants of 3-manifolds
associated with classical simple Lie algebras}, Int. J Math.
\textbf{4} (1993), 323-358.
\bibitem[TW2]{TuraevWenzl2} V.G. Turaev and H. Wenzl, \emph{Semisimple and modular
tensor categories from link invariants}, Math. Ann. \textbf{309}
(1997) 411-461.
\bibitem[TuW1]{TuW1} I. Tuba and H. Wenzl,
\emph{Representations of the braid groups $B_3$ and of
$SL(2,\Z)$}, Pacific J. Math. \textbf{197} (2001), no. 2, 491-510.
\bibitem[TuW2]{TuW2} I. Tuba and H. Wenzl, \emph{On braided tensor categories of type BCD},
to appear, J. Reine Agnew. Math.
\bibitem[W1]{BCD} H. Wenzl, \emph{Quantum groups and subfactors of
Lie type B, C and D}, Comm. Math. Phys. \textbf{133} (1990)
383-433.
\bibitem[W2]{wenzl} H. Wenzl, \emph{$C^*$ tensor categories from quantum groups},
J of AMS, \textbf{11} (1998) 261-282.
\bibitem[W3]{Wenzl88} H. Wenzl, \emph{Hecke algebras of type $A_n$ and subfactors},
 Invent.
Math. \textbf{92} (1988) 349-383.
\bibitem[Wy]{Wy} H. Weyl, The Classical Groups; Their Invariants
and Representations, Princeton University Press,
Princeton, 1939.
\end{thebibliography}

\end{document}